\documentclass[12pt]{article}
 \usepackage{graphicx}
\newfam\Bbbfam
\font\tenBbb=msbm10 scaled 1200
\font\sevenBbb=msbm8
\font\fiveBbb=msbm6
\textfont\Bbbfam=\tenBbb
\scriptfont\Bbbfam=\sevenBbb
\scriptscriptfont\Bbbfam=\fiveBbb
\def\Bbb{\fam\Bbbfam \tenBbb}

\begin{document}%

\newcounter{thm}[section]
\renewcommand{\thethm}{\arabic{section}.\arabic{thm}}
\newtheorem{theo}[thm]{Theorem}
\newtheorem{pro}[thm]{Proposition}
\newtheorem{lem}[thm]{Lemma}
\newtheorem{cor}[thm]{Corollary}

\makeatletter
  \renewcommand{\theequation}{% 
       \thesection.\arabic{equation}}
  \@addtoreset{equation}{section}
\makeatother

\title{Hitting time of a half-line by a two-dimensional 
nonsymmetric random walk} 
\author{Yasunari Fukai} 

\date{ }
\maketitle%

\def\thefootnote{ }

\pagestyle{myheadings}
\markboth{}{ }

\footnote{Yasunari Fukai: Faculty of Mathematics, Kyushu University,
Nishi-Ku, Fukuoka 812-0395, JAPAN.}

\footnote{e-mail: fukai@math.kyushu-u.ac.jp}



\begin{abstract}
We consider the probability that a two-dimensional 
random walk starting from the origin never returns to the 
half-line $ ( - \infty ,0] \times \{ 0 \}$ before time $n$. 
Let $X^{(1)}=(X_{1},X_{2})$ be the increment of the 
two-dimensional 
random walk. For an aperiodic random walk with moment conditions 
($E[X_{2}]=0 $ and $ E[ \vert X_{1} \vert^{\delta} ]< \infty,
E[ \vert X_{2} \vert^{2+ \delta} ]< \infty $ for some $ \delta \in (0,1)$), 
we obtain an asymptotic estimate (as $n \rightarrow \infty $) 
of this probability by assuming the behavior of the characteristic function 
of $X_{1}$ near zero. 
\end{abstract}

\section{Introduction and Statement of Results}

Let us consider the probability that a two-dimensional 
random walk starting from the origin never returns to the 
half-line $ ( - \infty ,0] \times \{ 0 \}$ before time $n$. 
In this paper we investigate the asymptotic behavior (as 
$ n \rightarrow \infty $ ) of this probability.

For a simple random walk, this probability, multiplied by $ n^{1/4}$, 
is bounded both from above and from below by positive constants. 
This bound was shown by Lawler [5,(2.35)], and  
the estimate plays an important role in 
computing an upper bound of a hitting distribution 
from infinity for a simple random walk 
(see [5,Proposition 2.4.10]). The upper bound 
gives a growth estimate for the 
two-dimensional diffusion limited aggregation 
(DLA in short) model (see [4]) and an estimate of the 
two-dimensional 
intersection exponent (see [6]).

Several authors have already studied 
the same probability that we consider here.
In [1], Bousquet-M$\acute{\mbox{e}}$lou and Schaeffer 
treated two-dimensional random walks  
satisfying the condition that the increment of the random 
walk takes finite values with equal probability by using a combinatorial 
argument. In [3], Isozaki treated two-dimensional symmetric 
random walks. 
In [2], a class of two-dimensional aperiodic random walks which satisfy good 
conditions (a mean of zero and the 
$2+ \delta (>2) $-th absolute moment exists) 
such as the simple random walk is considered. 
The present approach is probabilistic. We mainly treat 
two-dimensional nonsymmetric random walks.

Let $ \{ X^{(k)} \}_{k=1}^{ \infty }$ be a sequence of 
independent identically distributed random variables 
with values in ${\Bbb Z}^{2}$. A two-dimensional random 
walk $\{S(n)\}_{n=0}^{ \infty }$ starting at $x \in 
{\Bbb Z}^{2}$ is defined by 
$$ S(0)=x \quad and \quad S(n)=x + \sum_{k=1}^{n} X^{(k)}.$$
If the random variable $X^{(1)}$ takes the four values $(1,0),(-1,0),(0,1),
(0,-1)$ with equal probability, the random walk is called simple. We write 
$X^{(1)} =(X_{1},X_{2})$ and $S(n)$ $= (S_{1}(n),S_{2}(n))$. 
We denote by $P_{x}$ the probability law of the two-dimensional random 
walk starting at $x$ and by $E_{x}$ the expectation with respect to $P_{x}$.
For a subset $A$ of ${\Bbb Z}^{2}$, define 
$$ \tau _{A}= \inf \{ n \geq 1 : S(n) \in A \},$$
(note that the time $n=0$ is not included here). 
We also denote by $U$ the first coordinate axis:
$$ U= \{ x =(x_{1},x_{2}) \in {\Bbb Z}^{2} : x_{2}=0 \} \quad \quad \quad $$
and by $V_{-}$ and $V_{+}$ the left half and right half of $U$, respectively:
$$ V_{-}= \{ x =(x_{1},x_{2}) \in {\Bbb Z}^{2} : x_{1} \leq 0 , x_{2}=0 \} ,$$
$$ V_{+}= \{ x =(x_{1},x_{2}) \in {\Bbb Z}^{2} : x_{1} \geq 0 , x_{2}=0 \} .$$

\begin{theo}
Assume that a two-dimensional random walk satisfies the 
following conditions: 
$$(a) \ \mbox{aperiodic }( i.e.,\mbox{ the smallest additive subgroup 
containing } $$
$$ \{ x \in {\Bbb Z}^{2} \ | \ P \{X^{(1)} = x \} >0 \} 
\mbox{ agrees with } {\Bbb Z}^{2}), \quad \quad \quad \quad $$
$$(b) \ E[X_{2}]=0, \quad \quad \quad \quad \quad \quad \quad \quad
 \quad \quad \quad \quad \quad \quad \quad \quad \quad \quad \quad \quad $$
$$(c) \ E[ \vert X_{1} \vert ^{ \delta } ] < \infty , \quad 
E[ \vert X_{2} \vert ^{2 + \delta } ] < \infty \mbox{ for some }
\delta \in (0,1). \quad \quad $$
If there exist constants 
$ \alpha \in (0,2], \varepsilon >0, c_{1}>0,c_{2}>0$ such that 
$$ 1-E[ \cos \theta_{1}X_{1}] - \frac{(E[( \sin \theta_{1}X_{1})X_{2} ] )^{2}}
{2E[X_{2}^{2}]} 
 + \frac{(E[( -1 + \cos \theta_{1}X_{1})X_{2} ] )^{2}}
{2E[X_{2}^{2}]} $$
$$ = c_{1} \vert \theta_{1} \vert ^{ \alpha } + 
O ( \vert \theta_{1} \vert ^{ \alpha + \varepsilon } ) \quad 
( \theta_{1} \rightarrow 0)  \quad \quad \quad \quad \quad 
\quad \quad \quad \quad \quad \quad $$
and 
$$ E[ \sin \theta_{1}X_{1}] - \frac{ E[ (\sin \theta_{1}X_{1})X_{2}] E[(-1+ \cos\theta_{1}X_{1} )X_{2}]}{ E[X_{2}^{2}]} \quad \quad \quad \quad $$
$$ = c_{2} \vert \theta_{1} \vert^{ \alpha } 
\frac{ \theta_{1}}{ \vert \theta_{1} \vert } + O ( \vert \theta_{1} \vert ^{ \alpha + \varepsilon } )  \quad ( \theta_{1} \rightarrow 0) \quad 
\quad \quad \quad \quad \quad \quad \quad $$
then 
\begin{equation}
  P_{0} \{ \tau_{V_{-}} >n \} \asymp n^{-1/4 + \beta } 
\label{eq:200}
\end{equation}
where 
$$ \beta = \frac{1}{ \alpha \pi } \arcsin  \frac{ c_{2}}
{ \sqrt{2} \sqrt{c_{1}^{2}+c_{2}^{2}+ c_{1} \sqrt{c_{1}^{2}+c_{2}^{2}}}} $$
and we write $ a_{n} \asymp b_{n} $ if there exist positive constants 
$ \tilde{c}_{*},\tilde{c}^{*}$ such that 
$$ \tilde{c}_{*} b_{n} \leq a_{n} \leq \tilde{c}^{*} b_{n}.$$ 
\end{theo}

As an example, we consider the increment $X^{(1)}=(X_{1},X_{2})$ 
of a two-dimensional random walk which satisfies the following 
conditions: 
$$ \mbox{ (I) aperiodic, } \quad \quad \quad \quad \quad \quad
 \quad \quad \quad \quad \quad \quad \quad \quad \quad \quad
 \quad \quad \quad \quad \quad \quad \quad \quad \quad \quad \quad
  \quad $$
$$ \mbox{(II) } X_{1} \mbox{ and } X_{2} \mbox{ are independent, } \quad \quad \quad 
 \quad \quad \quad \quad \quad \quad \quad \quad \quad \quad \quad 
 \quad \quad \quad \quad \quad $$
$$ \mbox{(III) } E[X_{2}]=0 \mbox{ and } E[ \vert X_{2} \vert ^{2+ \delta }] < \infty 
 \mbox{ for some } \delta >0,  \quad \quad \quad \quad \quad \quad 
 \quad \quad \quad \quad \quad $$
$ \quad \quad \quad \mbox{(IV) }$ for $ \alpha \in (1,2)$, 
\[ 
P \{ X_{1}=n \} = \left\{ 
\begin{array}{ll} 
c_{-} \frac{1}{ \vert n \vert ^{1+ \alpha }} & ( n \in \{ -1,-2 , \cdots \}) \\
 c_{+} & (n=1) \\
0 & ( \mbox{ otherwise } ), \\
\end{array} 
\right. 
\]
where the constants $c_{-} $ and $c_{+}$ are assumed to be adjusted so that 
$$ 
\sum_{ n \in {\Bbb Z}} P \{ X_{1}=n\} =1 \  \mbox{ and } \ 
 \sum_{ n \in {\Bbb Z}} n P \{ X_{1}=n \} =0 $$ 
hold. By simple calculus, 
$$ 1- E[ \cos \theta_{1}X_{1}] = \frac{ c_{-} \pi  }{ 2 \Gamma ( \alpha +1) \cos \left( ( \alpha -1) \pi /2 \right) }  \vert \theta_{1} \vert ^{ \alpha} 
+ O( \vert \theta_{1} \vert ^{2} ) ,$$
$$ E[ \sin \theta_{1}X_{1}] = \frac{ c_{-} \pi  }{ 2 \Gamma ( \alpha +1) \sin \left( ( \alpha -1) \pi /2 \right) } 
\vert \theta_{1} \vert ^{ \alpha} 
\frac{  \theta_{1}  }{ \vert \theta_{1} \vert } 
+ O( \vert \theta_{1} \vert ^{2} ) ,
$$ 
where $ \Gamma $ is the gamma function. Thus 
$$ \beta= \frac{2 - \alpha }{ 4 \alpha} $$
and 
$$ P_{0} \{ \tau_{V_{-}} >n \} \asymp n^{-( \alpha -1)/(2 \alpha) } . $$

\begin{cor}
Assume $(a)$, $(b)$ and $(c)$. If there exist constants 
$ \alpha \in (0,2], \varepsilon >0, c_{1}>0 $ such that 
$$ 1-E[ \cos \theta_{1}X_{1}] - \frac{(E[( \sin \theta_{1}X_{1})X_{2} ] )^{2}}
{2E[X_{2}^{2}]}  + \frac{(E[( -1 + \cos \theta_{1}X_{1})X_{2} ] )^{2}}
{2E[X_{2}^{2}]} $$
$$ = c_{1} \vert \theta_{1} \vert ^{ \alpha } + 
O ( \vert \theta_{1} \vert ^{ \alpha + \varepsilon } ) 
 \quad ( \theta_{1} \rightarrow 0) \quad 
\quad \quad \quad \quad \quad \quad \quad \quad \quad \quad $$
and 
$$ E[ \sin \theta_{1}X_{1}] - \frac{ E[ (\sin \theta_{1}X_{1})X_{2}] E[(-1+ \cos\theta_{1}X_{1} )X_{2}]}{ E[X_{2}^{2}]} \quad \quad \quad \quad \quad \quad $$
$$ = O ( \vert \theta_{1} \vert ^{ \alpha + \varepsilon } )  \quad ( \theta_{1} \rightarrow 0) \quad \quad \quad \quad \quad \quad \quad \quad
 \quad \quad \quad \quad \quad \quad \quad   $$
then 
\begin{equation}  
P_{0} \{ \tau_{V_{-}} >n \} \asymp n^{-1/4 } .
\label{eq:210}
\end{equation}
\end{cor}

For a two-dimensional simple random walk, (\ref{eq:210}) was shown by 
Lawler [5,(2.35)] and Kesten [4]. Bousquet-M$ \acute{\mbox{e}}$lou and 
Schaeffer [1] considered $P_{0} \{ \tau_{V_{-}} >n \}$ for a two-dimensional 
random walk satisfying the condition that the increment of the random walk 
takes finite values with equal probability. For a two-dimensional simple 
random walk, they computed 
$$ \lim_{n \rightarrow \infty } P_{0} \{ \tau_{V_{-}} >n \}/ 
\left( \frac{ \sqrt{ 1+ \sqrt{2}}}{ 2 \Gamma (3/4)} n^{-1/4} \right) =1 $$
by using a combinatorial argument. 
We consider an aperiodic random walk 
with a mean of zero and a finite $2+ \delta (>2)$-th absolute moment. 
In [2], it is proved that for this random walk, 
$ n^{1/4}  P_{0} \{ \tau_{V_{-}} >n \}$ converges to some positive constant 
as $ n \rightarrow \infty $. 
This random walk satisfies the assumptions of Corollary 1.2 and we have 
$ P_{0} \{ \tau_{V_{-}} >n \} \asymp n^{-1/4 } $. This result is weaker than 
that in [2], but Corollary 1.2 is applicable to a large class of random walks.

\begin{cor}
Assume $(a)$, $(b)$ and $(c)$. 
If there exist constants 
$ \alpha \in [1,2], \varepsilon >0, c_{2}>0$ such that 
$$ 1-E[ \cos \theta_{1}X_{1}] - \frac{(E[( \sin \theta_{1}X_{1})X_{2} ] )^{2}}
{2E[X_{2}^{2}]} 
 + \frac{(E[( -1 + \cos \theta_{1}X_{1})X_{2} ] )^{2}}
{2E[X_{2}^{2}]} $$
$$ = O ( \vert \theta_{1} \vert ^{ \alpha + \varepsilon } ) \quad 
( \theta_{1} \rightarrow 0) \quad 
\quad \quad \quad \quad \quad \quad \quad \quad \quad \quad 
 \quad \quad \quad \quad $$
and 
$$ E[ \sin \theta_{1}X_{1}] - \frac{ E[ (\sin \theta_{1}X_{1})X_{2}] E[(-1+ \cos\theta_{1}X_{1} )X_{2}]}{ E[X_{2}^{2}]}  \quad \quad \quad \quad 
 \quad \quad $$
$$ = c_{2} \vert \theta_{1} \vert^{ \alpha } 
\frac{ \theta_{1}}{ \vert \theta_{1} \vert } + O ( \vert \theta_{1} \vert ^{ \alpha + \varepsilon } ) \quad ( \theta_{1} \rightarrow 0) 
\quad \quad \quad \quad \quad 
\quad \quad \quad \quad $$
then 
$$  P_{0} \{ \tau_{V_{-}} >n \} \asymp n^{-1/4 + 1/(4 \alpha ) } . $$
\end{cor}

We consider the increment $X^{(1)}=(X_{1},X_{2})$ 
of a two-dimensional random walk which satisfies the following 
conditions: (I) aperiodic, (II) $X_{1} $ and $X_{2}$ are independent, 
(III) $E[X_{2}]=0$, and 
$E[ \vert X_{2} \vert ^{2+ \delta }] < \infty $ for some $ \delta >0$, 
(IV) '
$ P \{ X_{1}= 1 \}=1. $ It is clear that for $ n \in {\Bbb N}$,
$$ P_{0} \{ \tau_{V_{-}} >n \} = 1 .$$
We apply Corollary 1.3 to this random walk. By simple calculus,  
$$ 1-E[ \cos \theta_{1} X_{1}]= O( \vert \theta_{1} \vert ^{2} ), \quad 
E[ \sin \theta_{1} X_{1}]= \theta_{1} + O( \vert \theta_{1} \vert ^{2} ) $$
and thus 
$$ P_{0} \{ \tau_{V_{-}} >n \} \asymp 1, $$
which is in agreement with $ P_{0} \{ \tau_{V_{-}} >n \} = 1 .$

We consider $P_{0} \{ T_{ \lambda } < \tau_{V_{-}} \} $, 
where $T_{ \lambda } $ is a random variable independent of the random walk with
the law 
$P \{ T_{ \lambda } =j
\}= (1 - \lambda ) \lambda ^{j} \  (j=0,1,2,\cdot \cdot \cdot )$. By the
Tauberian theorem, the 
asymptotic behavior of $P_{0} \{ T_{ \lambda }  < \tau_{V_{-}} \}$ as 
$ \lambda \uparrow 1$ 
determines the asymptotic behavior of $P_{0} \{ \tau_{V_{-}} > n \}$ as $ n
\rightarrow \infty$.

If we can calculate the 
asymptotic behavior of $ P_{0} \{ T_{ \lambda }  < \tau_{V_{-}} \} \times P_{0} \{ T_{ \lambda }  < \tau_{V_{+}} \} $ and $ P_{0} \{ T_{ \lambda }  < \tau_{V_{-}} \} / P_{0} \{ T_{ \lambda }  < \tau_{V_{+}} \} $ 
as $ \lambda \uparrow 1 $, then we have the 
asymptotic behavior of $P_{0} \{ T_{ \lambda } < \tau_{V_{-}} \} $ 
as $ \lambda \uparrow 1 $. This method is based on [2]. 
However, to have Theorem 1.1, we need to obtain an estimate of 
$ P_{0} \{ T_{ \lambda }  < \tau_{V_{-}} \} / P_{0} \{ T_{ \lambda }  < \tau_{V_{+}} \} $ which is different from that in [2].

In Section 2, we state a certain relation between $ P_{0} \{ T_{ \lambda }  < \tau_{V_{-}} \} \times P_{0} \{ T_{ \lambda }  < \tau_{V_{+}} \} $ and 
$ P_{0} \{ T_{ \lambda }  < \tau_{U} \} $ and give the 
asymptotic behavior of $ P_{0} \{ T_{ \lambda }  < \tau_{U} \} $ as $ \lambda 
\uparrow 1.$ We show that $ P_{0} \{ T_{ \lambda }  < \tau_{V_{-}} \} / 
P_{0} \{ T_{ \lambda }  < \tau_{V_{+}} \} $ may be expressed using 
the characteristic function of the increment of the random walk.

In Section 3, we state some results needed to prove Theorem 1.1.

In Section 4, we give a proof of Theorem 1.1.

In the appendix, we give the proof of lemmas not proved in Section 3.

\section{Fundamental relations}

\begin{lem}
For an arbitrary two-dimensional random walk, 
$$ P_{0} \{ T_{ \lambda } < \tau _{U} \} = 
P_{0} \{ T_{ \lambda } < \tau _{\tilde{V}_{+} }
\}
 P_{0} \{ T_{ \lambda } < \tau _{V_{-} } \}  \quad \quad \quad \quad  \quad
\quad 
 \quad \quad \quad  \quad \quad \quad  \quad $$
$$ \quad \quad \quad \quad \quad \quad \quad
 = ( 1- P_{0} \{ S( \tau _{V_{+}} ) = 0 , \tau_{V_{+}} \leq T_{ \lambda } \}
)^{-1} P_{0} \{ T_{ \lambda } < \tau _{V_{+} } \}  P_{0} \{ T_{ \lambda } <
\tau _{V_{-} }
\} , $$ where $ \tilde{V}_{+} = V_{+} \setminus \{0 \}$.
\end{lem}

The proof of Lemma 2.1 is identical to that given for Proposition 2.4.6 
in [5] and is thus omitted (see [2, Lemma 2.1]).

\begin{lem}
A two-dimensional random walk satisfying (a), (b) and $E[X_{2}^{2}] < 
\infty $ has the property that 
$$ P_{0} \{ T_{ \lambda } < \tau _{U} \} \sim ( 2E[X_{2}^{2}] )^{1/2} 
(1  - \lambda )^{1/2} ,$$
where we write $ \alpha_{ \lambda } \sim \beta_{ \lambda } $ if 
$ \lim_{ \lambda \uparrow 1 } \alpha_{ \lambda }/ \beta_{ \lambda }=1.$ 
\end{lem}

The proof of Lemma 2.2 is identical to that given 
in [5, p.68] and is thus omitted (see [2, Proposition 3.1]).

To formulate the next lemma, we introduce some notation. For a two-dimensional 
random walk satisfying $(b)$ and $E[X_{2}^{2}] < \infty $, we introduce the 
following random variables: 
$$ \eta(1)= \inf \{ j \geq 1 \ | \ S(j) \in U \}, \quad 
\zeta(1)= S_{1}(\eta(1)), $$
$$ \eta(k+1)= \inf \{ j > \eta(k) \ | \ S(j) \in U \} \quad 
\quad \quad \quad \quad $$
$$ \zeta(k+1) = S_{1}( \eta(k+1)), \quad (k= 1,2, \cdot \cdot \cdot ). 
\quad \quad $$
For  $ \lambda \in (0,1)$, we define 
$$ c( \lambda ) = \exp \left( - \sum_{k=1}^{ \infty }
 \frac{1}{k} E_{0}[ \lambda ^{ \eta(k)} ; \zeta(k)=0 ] \right), $$
$$ f_{ \infty } (z; \lambda ) = \exp \left( - \sum_{k=1}^{ \infty }
 \frac{1}{k} E_{0}[ \lambda ^{ \eta(k) } z ^{ \zeta(k)} ; \zeta(k) <0 ] \right), $$
$$ f_{0} (z; \lambda ) = \exp \left( - \sum_{k=1}^{ \infty }
 \frac{1}{k} E_{0}[ \lambda ^{ \eta(k)} z ^{ \zeta(k)} ; \zeta(k) >0 ] \right). $$

\begin{lem}
For every two-dimensional random walk satisfying $(b)$ and 
$E[X_{2}^{2}] < \infty $, 
$$ P_{0} \{ T_{ \lambda } < \tau_{V_{-}} \}= c ( \lambda )f_{ \infty }(1; \lambda ) ,$$
$$ P_{0} \{ T_{ \lambda } < \tau_{V_{+}} \}= c ( \lambda )f_{ 0 }(1; \lambda ).
$$
Moreover, for 
$ l \in {\Bbb Z} $,
\begin{equation}
 E_{0}[ \lambda ^{\eta(k)}; \zeta(k)=l ] = \frac{1}{ 2 \pi } 
\int_{ - \pi }^{ \pi } e^{- i \theta_{1} l } 
\left( 1- \left( \frac{1}{2 \pi} \int_{ -\pi}^{\pi} \frac{1}
{ 1 - \lambda \phi ( \theta )} \ d \theta_{2} \right) ^{-1} \right)^{k}
\ d \theta_{1} . \label{eq:2.1}
\end{equation}
\end{lem}

The proof of Lemma 2.3 is identical to that given in \S 17 Proposition 5 
of [7] and is thus omitted (see [2, Proposition 2.4]).

\section{Preliminary lemmas}

In this section, we will consider $a_{ \lambda }( \theta_{1} )$ and 
$b_{ \lambda } ( \theta_{1} )$ defined by 
$$ a_{ \lambda }( \theta_{1} )= \Re \left[ \frac{1}{2 \pi } 
\int_{ - \pi }^{ \pi } \frac{1}{ 1 - \lambda \phi ( \theta_{1} ,
\theta_{2}) } \ d \theta_{2} \right] $$
and 
$$ b_{ \lambda }( \theta_{1} )= \Im \left[ \frac{1}{2 \pi } 
\int_{ - \pi }^{ \pi } \frac{1}{ 1 - \lambda \phi ( \theta_{1} ,
\theta_{2}) } \ d \theta_{2} \right]. $$
Let
$$ \varphi _{ \lambda } ( \theta_{1}, \theta_{2} ) = 
1 - \lambda E[ \cos \theta_{1} X_{1}] + \lambda E[( \sin \theta_{1} 
X_{1})X_{2} ] \theta_{2} + \frac{ \lambda }{2} E[X_{2}^{2} ] 
\theta_{2}^{2} $$
$$ \quad \quad  -i \{ \lambda E[ \sin \theta_{1} X_{1}] + \lambda 
E[( -1+ \cos \theta_{1} X_{1})X_{2}] \theta_{2} \}. $$
Then we have the following lemma.

\begin{lem} 
Assume $(a)$, $(b)$ and $(c)$ hold. Then 
there exists a constant $c_{3}>0$ such that, for $ \lambda \in (1/2,1)$ 
and $ \theta =( \theta_{1}, \theta_{2} ) \in [ - \pi ,\pi ]
\times [ - \pi ,\pi ], $
\begin{equation} 
 \left\vert \frac{ \varphi _{ \lambda } (  \theta ) }
{ 1 - \lambda \phi ( \theta ) } - 1 \right\vert \leq c_{3} ( 
\vert \theta_{1} \vert ^{ \hat{\delta}} + 
\vert \theta_{2} \vert ^{ \delta } ) ,
\label{eq:vp} 
\end{equation}
where $ \hat{ \delta }= \delta ^{2}/(2 + \delta ).$

In addition, there exist $ \tilde{c}_{3} >0 $ and $ r_{0} >0$ 
such that, for $\lambda \in (1/2,1)$ 
and $ \theta =( \theta_{1}, \theta_{2} ) \in $ $ [ - \pi ,\pi ]
\times [ - \pi ,\pi ] $ with $ \vert \theta \vert < r_{0} $, 
\begin{equation} 
 \left\vert \frac{ \Re[ 1- \lambda \phi ( \theta )] }
{ \Re [ \varphi _{ \lambda } (  \theta ) ]} -1 \right\vert 
\leq \tilde{c}_{3}  ( \vert \theta_{1} \vert ^{ \hat{\delta}} + 
\vert \theta_{2} \vert ^{ \delta } ) ,
\label{eq:rvp} 
\end{equation}
where $ \hat{ \delta } $ is the same constant as in (\ref{eq:vp}).
\end{lem}

$Proof.$ By the definition of the characteristic function, 
$$ 1 - \lambda \phi ( \theta ) = 1 - \lambda E[e^{i \theta_{1}X_{1}} 
e^{i \theta_{2} X_{2}}]. $$
The inequality 
$$ \left\vert e^{ix} - \left( 1+ix- \frac{1}{2}x^{2} \right) \right\vert 
\leq 2 \vert x \vert ^{2+ \delta } \quad ( x\in {\Bbb R}) $$
and the condition $E[\vert X_{2} \vert ^{2 + \delta }] < \infty $ give 
$$ \left\vert E[e^{i \theta_{1}X_{1}} 
e^{i \theta_{2} X_{2}}] - E \left[ e^{i \theta_{1} X_{1}} \left( 
1+i \theta_{2} X_{2} - \frac{1}{2} \theta_{2}^{2} X_{2}^{2} \right) 
 \right] \right\vert \leq 2 E[ \vert X_{2} \vert^{2+ \delta }  ] 
\vert \theta_{2} \vert ^{2 + \delta }. $$
Using the H$ \ddot{ \mbox{o}} $lder inequality and condition 
$(c)$, $ E[ \vert X_{1} \vert ^{ \hat{ \delta }}
 X_{2}^{2} ] < \infty .$ The inequality $ \vert e^{ix} -1  \vert $ 
$ \leq 3 \vert x \vert ^{ \hat{ \delta } } \ (x \in {\Bbb R})$ gives 
$$ \vert E[ e^{i \theta_{1}X_{1}} X_{2}^{2} ] - E[ X_{2}^{2}] 
\vert \leq 3 E[\vert X_{1} \vert ^{ \hat{ \delta }} X_{2}^{2} ] 
\vert \theta_{1} \vert^{ \hat { \delta }} .$$
We notice that 
$$ E[( -1+ \cos \theta_{1} X_{1})X_{2}] = E[
( \cos \theta_{1} X_{1})X_{2}] $$
from condition (b).
Hence 
\begin{equation}
 \vert 1 - \lambda \phi ( \theta ) - \varphi _{ \lambda } ( \theta )
\vert \leq c ( \vert \theta _{1} \vert ^{ \hat{ \delta }} 
 + \vert \theta _{2} \vert ^{\delta } ) 
\theta_{2}^{2} 
\label{eq:1} 
\end{equation}
for some suitable constant $c>0$. Since the random walk is aperiodic, 
there exists a constant $c_{*}>0$ such that, for 
$ \lambda \in (1/2,1)$ and $ \theta \in [- \pi ,\pi ] \times 
[- \pi ,\pi ], $
\begin{equation}
\vert 1 - \lambda \phi ( \theta ) \vert \geq  \Re [ 1 - 
\lambda  \phi ( \theta )] \geq 
\lambda (1 - \Re [ \phi ( \theta )] ) \geq \frac{1}{2}c_{*} 
\vert \theta \vert ^{2} 
\label{eq:2} 
\end{equation}
as shown in [7, \S7 Proposition 5]. (\ref{eq:vp}) follows from 
(\ref{eq:1}) and (\ref{eq:2}).

(\ref{eq:1}) implies that 
\begin{equation}
 \vert \Re[ 1- \lambda \phi ( \theta ) ] - \Re[   \varphi _{ \lambda } ( \theta ) ] \vert \leq c ( \vert \theta _{1} \vert ^{ \hat{ \delta }} 
 + \vert \theta _{2} \vert ^{\delta } ) 
 \theta_{2}^{2} .
\label{eq:11} 
\end{equation}

By combining (\ref{eq:11}) with the last inequality in (\ref{eq:2}), 
there exists $ r_{0}>0 $ (small enough) such that, for 
$ \lambda $ $ \in (1/2,1)$ 
and $ \theta =( \theta_{1}, \theta_{2} ) \in $ $ [ - \pi ,\pi ]
\times [ - \pi ,\pi ] $ with $ \vert \theta \vert < r_{0} $, 
\begin{equation}
 \Re[ \varphi _{ \lambda } ( \theta ) ]  \geq \frac{1}{4}c_{*} 
\theta_{2}^{2} .
\label{eq:22} 
\end{equation}
(\ref{eq:rvp}) follows from 
(\ref{eq:11}) and (\ref{eq:22}).

To estimate $a_{ \lambda }( \theta )$, 
we calculate 
$$ \Re \left[ \frac{1}{ 2 \pi} \int_{ - r_{0} /2 }^{ r_{0}/2 } \frac{1}{ 
\varphi _{ \lambda } (  \theta_{1}, \theta_{2} )  } \ d \theta_{2}
\right]  = \frac{1}{ 2 \pi} \int_{ - r_{0} /2 }^{ r_{0}/2 } \Re \left[
 \frac{1}{ 
\varphi _{ \lambda } (  \theta_{1}, \theta_{2} )  } \right] \ d \theta_{2} , $$
where $r_{0}$ is the same constant as in Lemma 3.1. 
Let 
$$ A= \frac{ E[ ( \sin \theta_{1}X_{1}) X_{2}] }{ (1/2) E[X_{2}^{2} ] } ,
\quad B= \frac{ 1 - \lambda E[ \cos \theta_{1}X_{1}] }
{( \lambda /2) E[X_{2}^{2} ] } ,$$
$$ C= \frac{ E[ ( -1 + \cos  \theta_{1}X_{1} )X_{2} ]}{ (1/2) E[X_{2}^{2} ] } ,
\quad D = \frac{ E[ \sin \theta_{1}X_{1} ] }{  (1/2) E[X_{2}^{2} ] }. $$
Then 
$$  \varphi _{ \lambda } 
(  \theta_{1}, \theta_{2} )  =
 \frac{ \lambda }{2} E[X_{2}^{2} ] \{ \theta_{2}^{2} + A \theta _{2} +B
 - i (C \theta_{2}+D ) \}  $$
and 
$$  \Re \left[
 \frac{1}{ 
\varphi _{ \lambda } (  \theta_{1}, \theta_{2} )  } \right] = 
\frac{ \Re[ \varphi _{ \lambda } (  \theta_{1}, \theta_{2} ) ]}
{ \vert \varphi _{ \lambda } (  \theta_{1}, \theta_{2} ) \vert^{2}} \ 
\mbox{ is a rational function.}$$
Put 
$$ K= ( 4B -A^{2}+C^{2})^{2}+ 4 (2D-AC)^{2} .$$
To calculate the integral of $ \Re [ 1/ \varphi _{ \lambda } (  \theta_{1}, \theta_{2} ) ] $ over $ [-r_{0}/2, r_{0}/2] $ with respect to $ \theta_{2}$, we use routine integration of a rational function using partial fractions. 
This calculation gives the following lemma.

\begin{lem}
Suppose a two-dimensional random walk satisfies  
$ 0 < E[X_{2}^{2} ] < \infty $. Then there exist 
$ c_{4}>0, \lambda_{*} \in (1/2 ,1) $ and $ s_{*}>0 $ such that, 
for all $ \lambda \in ( \lambda_{*},1) $ and $ \theta_{1} \in [ - \pi ,\pi ]$ 
with $ 0 < \vert \theta_{1} \vert < s_{*} $,  
\begin{equation}
 \left\vert  \frac{1}{ 2 \pi} \int_{ - r_{0} /2 }^{ r_{0} /2 } \Re \left[ \frac{1}{ 
\varphi _{ \lambda } (  \theta_{1}, \theta_{2} )  } \right] \ d \theta_{2} 
 - \tilde{a}_{ \lambda }( \theta_{1}) \right\vert \leq c_{4}, 
\label{eq:l2}
\end{equation}
where 
$$ \tilde{a}_{ \lambda }( \theta_{1}) =
\frac{ \sqrt{2}}{ \lambda E[X_{2}^{2}] \sqrt{K}} \sqrt{ 4B-A^{2}+C^{2} 
+ \sqrt{K}} . $$
\end{lem}

With the help of (\ref{eq:vp}),(\ref{eq:1}),(\ref{eq:2}) and (\ref{eq:22}), 
we have the following lemma. 
\begin{lem}
Suppose a two-dimensional random walk satisfies $(a)$, $(b)$ and $(c)$. Then 
there exist $ c_{5}>0 ,  \lambda _{*} 
\in (1/2,1)$ and $ s_{*} >0$ such that, for 
$ \lambda \in (  \lambda _{*} ,1)$ and $ \theta_{1} \in 
[ - \pi, \pi ] $ with $ 0< \vert \theta_{1} \vert < 
s_{*} $, 
$$
\left\vert a_{ \lambda }( \theta_{1}) - \frac{1}{2 \pi } 
\int_{-r_{0}/2}^{r_{0}/2}
\Re \left[\frac{1}{ \varphi_{ \lambda }( \theta_{1}, \theta_{2}) } \right]
\ d \theta_{2} \right\vert \leq c_{5}
 ( \vert \theta_{1} \vert ^{ \hat{\delta} } 
+ K^{ \delta /4 } \vert \log K \vert ) K^{ -1/4} .$$
\end{lem}

\begin{lem}
Suppose a two-dimensional random walk satisfies  
$ 0 < E[X_{2}^{2} ] < \infty $. Then there exist 
$ c_{6}>0, \lambda_{*} \in (1/2,1) $ and $ s_{*}>0 $ such that, 
for all $ \lambda \in ( \lambda_{*},1) $ and $ \theta_{1} \in ( - \pi ,\pi )$ 
with $ 0 < \vert \theta_{1} \vert < s_{*} $, 
\begin{equation}
 \left\vert  \frac{1}{ 2 \pi} \int_{ - r_{0} /2 }^{ r_{0} /2 } 
\Im \left[ \frac{1}{ 
\varphi _{ \lambda } (  \theta_{1}, \theta_{2} )  } \right] \ d \theta_{2}
 - \tilde{b}_{ \lambda }( \theta_{1} )
 \right\vert \leq c_{6} \left( \frac{ \vert 2D -AC \vert }{ \sqrt{K}}
+ \frac{ \vert C \vert }{K^{1/4}} \right),
\label{eq:33}
\end{equation}
where 
$$ \tilde{b}_{ \lambda }( \theta_{1} )= 
\frac{ \sqrt{2}}{ \lambda E[X_{2}^{2}] \sqrt{K}} 
\frac{ 2 (2D-AC) }{\sqrt{ 4B-A^{2}+C^{2} 
+ \sqrt{K}} } .$$
\end{lem}

\begin{lem} 
Suppose a two-dimensional random walk satisfies $(a)$, $(b)$ and $(c)$. Then 
there exist $ c_{7}>0 ,  \lambda _{*} 
\in (1/2,1)$ and $ s_{*} >0$ such that, for 
$ \lambda \in (  \lambda_{*} ,1)$ and $ \theta_{1} \in 
[ - \pi, \pi ] $ with $ 0< \vert \theta_{1} \vert < 
s_{*} $, 
$$ \left\vert b_{ \lambda } ( \theta_{1} ) - 
\frac{1}{2 \pi } \int_{-r_{0}/2}^{r_{0}/2} 
\Im \left[ \frac{1}{ \varphi _{ \lambda } (  \theta_{1}, \theta_{2} )  } \right] \ d \theta_{2} \right\vert  \quad \quad \quad \quad \quad \quad \quad
 \quad \quad \quad \quad \quad \quad \quad $$
$$  \quad \quad \quad 
\leq  c_{7} \{  \vert \theta_{1} \vert ^{ \hat{\delta}}+ 
\vert A \vert K^{(-1+ \delta)/4} + 
( \vert \theta_{1} \vert ^{\hat{ \delta }} + K^{ \delta /4} ) 
( \vert C \vert + \vert D \vert K^{ -1/4} ) K^{-1/4} \}  K^{-1/4} .$$
\end{lem}

The proofs of the above four lemmas are given in the appendix.

\section{Proof of Theorem 1.1}

We consider $ P_{0} \{ T_{ \lambda }  < \tau_{V_{-}} \} / P_{0} \{ T_{ \lambda }  < \tau_{V_{+}} \} $ under the assumption that the random walk satisfies 
(a),(b) and (c). By Lemma 2.3 and the definitions of 
$f_{0}(z; \lambda )$ and $f_{ \infty }(z; \lambda )$, this can be expressed as 
\begin{equation}
 \frac{ P_{0} \{ T_{ \lambda }  < \tau_{V_{-}} \} }
{ P_{0} \{ T_{ \lambda }  < \tau_{V_{+}} \}  } = 
\frac{ f_{ \infty }(1; \lambda )}{f_{0}(1; \lambda )} = 
 \exp \left( \lim_{L \rightarrow \infty } 
C_{L}( \lambda ) \right) ,
\label{eq:180}
\end{equation}
where 
$$ C_{L}( \lambda )=  \sum_{k=1}^{ \infty }
\frac{1}{k} \sum_{l=1}^{ L} \left( E_{0}[ \lambda ^{\eta(k)} ; \zeta(k)=l ] - 
E_{0}[ \lambda ^{\eta(k)} ; \zeta(k)=-l ] \right) . $$
Since the left-hand side of the first identity in (\ref{eq:180}) 
is a positive real 
number, 
\begin{equation}
 \exp \left( \lim_{L \rightarrow \infty } 
C_{L}( \lambda ) \right) = \exp \left( \lim_{L \rightarrow \infty } 
\Re[C_{L}( \lambda ) ] \right). 
\label{eq:190}
\end{equation}

To obtain an expression for $ \Re[C_{L}( \lambda )]$, 
we use (\ref{eq:2.1}) in Lemma 2.3 which implies that 
\begin{equation}
C_{L}( \lambda ) = \sum_{k=1}^{ \infty } \frac{1}{2 \pi } 
\int_{ - \pi }^{ \pi }  e_{L}( \theta_{1}) \frac{1}{k} \left( 1- \left( \frac{1}{2 \pi } 
\int_{ - \pi }^{ \pi } \frac{ 1}{ 1 - \lambda \phi ( \theta ) } \ 
d \theta_{2} \right)^{-1} \right)^{k} \ d \theta_{1} , 
\label{eq:clla}
\end{equation}
where 
$$ e_{L}( \theta_{1} )= \sum_{l=1}^{L} ( e^{-il \theta_{1}} - 
e^{i l \theta_{1}} ). $$
From the inequality  
$$ \sup_{ \theta_{1} \in [ - \pi ,\pi]} \left\vert 
1- \left( \frac{1}{ 2 \pi } \int_{ - \pi }^{ \pi }
\frac{ 1}{ 1 - \lambda \phi ( \theta ) } \ 
d \theta_{2} \right)^{-1} \right\vert <1,$$
which is shown in [2, p332-p333], we can interchange the order of summation 
and integration on the right-hand side of (\ref{eq:clla}) and then 
use the identity 
$$  - \sum_{k=1}^{ \infty } \frac{z^{k}}{k} = \mbox{Log} (1-z), \quad ( \vert z \vert <1). $$
Here, $ \mbox{Log} (z) $ is the principal logarithm of $z$ and 
we choose $ \vert \Im [ \mbox{Log} (z) ] \vert < \pi .$
$ C_{L} ( \lambda ) $ reduces to 
$$  \frac{1}{ 2 \pi } \int_{ - \pi }^{ \pi }
e_{L}( \theta_{1} ) \mbox{ Log } \left( 
\frac{1}{ 2 \pi } \int_{ - \pi }^{ \pi }
\frac{ 1}{ 1 - \lambda \phi ( \theta ) } \ 
d \theta_{2} \right) \ d \theta_{1} . $$

By a simple calculation, 
$$ e_{L}( \theta_{1}) = i \left\{ \frac{ \sin \theta_{1} }
{1 - \cos \theta_{1}} (-1+ \cos L \theta_{1}) - 
\sin L \theta_{1}  \right\} . $$
If $ \Re [z] >0 $, then $ \Im [ \mbox{Log} (z)] \in ( - \pi /2, \pi /2)$ and 
$$ \Im [ \mbox{Log} (z) ] = \arcsin \frac{ \Im [z]}{ \vert z \vert } . $$
Recall that $a_{ \lambda }( \theta_{1}) $ is the real part of $ (1/(2 \pi)) \int_{ -\pi }^{ \pi } 1/(1- \lambda \phi ( \theta )) d \theta_{2}$ and 
$ b_{ \lambda }( \theta_{1} )$ the imaginary part. It is easy to check that 
$ a_{ \lambda } ( \theta_{1} ) >0 $. 
Thus  $ \Re[ C_{L}( \lambda )] $ can be expressed as 
$$ \Re[ C_{L}( \lambda )] = \frac{1}{2 \pi } \int_{ - \pi }^{ \pi } 
\left\{ \frac{ \sin \theta_{1} }
{1 - \cos \theta_{1}} (1- \cos L \theta_{1}) + 
\sin L \theta_{1}  \right\} \quad \quad \quad \quad $$
$$  \quad \quad \quad \quad 
\times  \arcsin  
\frac{ b_{ \lambda }( \theta_{1})}{ \sqrt{ a_{ \lambda }^{2}( \theta_{1})+ 
b_{ \lambda }^{2}( \theta_{1})} }  \ d \theta_{1} . $$

Let 
$$ B ( \theta_{1})=  1-E[ \cos \theta_{1}X_{1}] - \frac{(E[( \sin \theta_{1}X_{1})X_{2} ] )^{2}}
{2E[X_{2}^{2}]} 
 + \frac{(E[( -1 + \cos \theta_{1}X_{1})X_{2} ] )^{2}}
{2E[X_{2}^{2}]} ,$$
$$ D( \theta_{1})= E[ \sin \theta_{1}X_{1}] - \frac{ E[ (\sin \theta_{1}X_{1})X_{2}] E[(-1+ \cos\theta_{1}X_{1} )X_{2}]}{ E[X_{2}^{2}]} . \quad \quad \quad
 \quad \quad \quad $$
Then $ \tilde{a}_{ \lambda }( \theta_{1}) $ 
in Lemma 3.2 can be written as 
$$ \tilde{a}_{\lambda }( \theta_{1}) = \frac{1}
{2 \sqrt{ \lambda E[X_{2}^{2}]}} \frac{1} { \sqrt{ (1- \lambda + \lambda 
B( \theta_{1}) )^{2}+ ( \lambda D( \theta_{1}) )^{2}}} \quad \quad 
\quad $$
$$ \quad \quad \quad \quad \quad 
\times \sqrt{ 1- \lambda + \lambda 
B( \theta_{1}) + \sqrt{(1- \lambda + \lambda 
B( \theta_{1}) )^{2}+ ( \lambda D( \theta_{1}) )^{2}}}, $$
and $ \tilde{b}_{ \lambda }
( \theta_{1} ) $ in Lemma 3.4 can be written as 
$$ \tilde{b}_{ \lambda }( \theta_{1}) = \frac{1}
{2 \sqrt{ \lambda E[X_{2}^{2}]}} \frac{ \lambda D( \theta_{1})} { \sqrt{ (1- \lambda + \lambda 
B( \theta_{1}) )^{2}+ ( \lambda D( \theta_{1}) )^{2}}} \quad \quad 
\quad $$
$$ 
\quad \quad \quad \quad \quad 
\times \frac{1}{ \sqrt{ 1- \lambda + \lambda 
B( \theta_{1}) + \sqrt{(1- \lambda + \lambda 
B( \theta_{1}) )^{2}+ ( \lambda D( \theta_{1}) )^{2}}}} . $$
Hence, 
$$ \frac{ \tilde{b}_{ \lambda }( \theta_{1}) }{ 
\sqrt{ \tilde{a}_{ \lambda }^{2}( \theta_{1})+ 
\tilde{b}_{ \lambda }^{2}( \theta_{1})}} \quad \quad
 \quad \quad \quad \quad \quad \quad \quad \quad \quad \quad
 \quad \quad \quad \quad \quad \quad \quad \quad \quad \quad 
 \quad \quad \quad \quad \quad \quad \quad \quad $$
$$ = \frac{ \lambda D( \theta_{1})}{ \sqrt{2} 
\sqrt{ (1- \lambda + \lambda 
B( \theta_{1}) )^{2}+ ( \lambda D( \theta_{1}) )^{2}+ 
(1- \lambda + \lambda 
B( \theta_{1}) ) \sqrt{(1- \lambda + \lambda 
B( \theta_{1}) )^{2}+ ( \lambda D( \theta_{1}) )^{2}}}}. $$

To prove Theorem 1.1, we recall the following assumption: 
\vspace{.2cm}

(A1) there exist constants $ \alpha 
\in (0,2] , \varepsilon >0, c_{1} >0, c_{2}>0 $ such that 
$$ B( \theta_{1}) = c_{1} \vert \theta_{1} \vert ^{ \alpha } + O ( 
\vert \theta_{1} \vert ^{ \alpha + \varepsilon } ) \quad 
( \theta_{1} \rightarrow 0) \quad \quad $$
and 
$$ D( \theta_{1} ) = c_{2} \vert \theta_{1} \vert ^{ \alpha }
\frac{ \theta_{1}}{ \vert \theta_{1} \vert } + 
O ( \vert \theta_{1} \vert ^{ \alpha + \varepsilon }) \quad ( \theta_{1} \rightarrow 0).$$
\vspace{.2cm}

From assumption (A1), it is easy to see that 
there exist constants $c_{8}>0,$ 
$s >0$ such that, for $ \lambda \in ( 1/2 ,1)$ and 
$ 0< \vert \theta_{1} \vert \leq s$, 
$$ \left\vert \sqrt{ (1- \lambda + \lambda 
B( \theta_{1}) )^{2}+ ( \lambda D( \theta_{1}) )^{2}} - 
 k_{ \lambda }( \theta_{1})^{1/2} \right\vert \quad 
 \quad \quad \quad \quad \quad \quad \quad $$
\begin{equation}
 \leq \left\vert B( \theta_{1}) - c_{1} \vert \theta_{1} \vert ^{ \alpha } 
\right\vert 
+ \left\vert D ( \theta_{1} ) - c_{2} \vert \theta_{1} \vert ^{ \alpha }
\frac{ \theta_{1}}{ \vert \theta_{1} \vert } \right\vert 
 \leq c_{8}  k_{ \lambda }( \theta_{1})^{1/2}
 \vert \theta_{1} \vert^{ \varepsilon } , 
\label{eq:bda}
\end{equation}
where 
$$ k_{ \lambda }( \theta_{1} )= (1 - \lambda + \lambda c_{1} \vert 
\theta_{1} \vert ^{ \alpha } )^{2}+ \left( \lambda c_{2} \vert \theta_{1} 
\vert ^{ \alpha } \frac{ \theta_{1}}{ \vert \theta_{1} \vert } \right)^{2}  .$$
(\ref{eq:bda}) implies the following lemma.
\begin{lem} 
Assume that $ 0< E[X_{2}^{2}] < \infty $ and $(A1)$ hold. Then 
there exist $ c_{9}>0, 
s_{*} >0$ such that, for $ \lambda \in ( 1/2 ,1)$ and 
$ 0< \vert \theta_{1} \vert \leq s_{*}$, 
$$ \left\vert \frac{ \tilde{b}_{ \lambda }( \theta_{1}) }{ 
\sqrt{ \tilde{a}_{ \lambda }^{2}( \theta_{1})+ 
\tilde{b}_{ \lambda }^{2}( \theta_{1})}} - 
\frac{ \lambda c_{2} \vert \theta_{1} \vert ^{ \alpha } ( \theta_{1} / 
\vert \theta_{1} \vert ) }{ \sqrt{2} 
Q_{ \lambda } ( \vert \theta_{1} \vert ^{ \alpha }) }  \right\vert \leq c_{9} \vert \theta_{1} \vert^{ \varepsilon } , $$
where 
$$ Q_{ \lambda }(s)= \left\{ (1- \lambda +\lambda c_{1}s )^{2}+ 
( \lambda c_{2}s )^{2} + 
( 1- \lambda +\lambda c_{1}s) \sqrt{ 
(1- \lambda +\lambda c_{1}s )^{2}+ ( \lambda c_{2}s )^{2} } 
\right\}^{1/2} .$$
\end{lem}
From (\ref{eq:22}) and the definitions of A and B, 
$$
\frac{1}{4}c_{*} \vert \theta_{2} \vert ^{2} \leq \Re [ \varphi_{ \lambda } 
( \theta_{1} , \theta_{2} ) ] = \frac{ \lambda }{2} E[X_{2}^{2}] 
\left\{ 
\left( \theta_{2} + \frac{A}{2} \right)^{2} + B - \frac{ A^{2}}{4} \right\} 
$$
for $ \lambda \in (1/2,1)$ 
and $ \vert (\theta_{1} ,\theta_{2} ) \vert \leq r_{0} .$ 
Since A tends to zero as $ \vert \theta_{1} \vert \rightarrow 0$, 
we choose $ \theta_{1} $ small enough so that 
$ \vert ( \theta_{1}, -A/2) \vert \leq r_{0} .$ 
Then the above inequality (setting $ \theta_{2}= -A/2$) gives 
$$ \frac{c_{*}}{2 \lambda E[X_{2}^{2}] }A^{2} \leq 4B -A^{2}, $$
and letting $ \lambda \uparrow 1$,  
$$ \frac{ c_{*}}{2 E[X_{2}^{2}] } A^{2} \leq 4 
\frac{1- E[ \cos \theta_{1}X_{1}]}{(1/2) E[X_{2}^{2}] } - A^{2} \leq 
8 \frac{B( \theta_{1})}{ E[X_{2}^{2}]} .$$
From $ 1 - E[ \cos \theta_{1} X_{1}] \geq \left( E[( \sin  \theta_{1} X_{1}) 
X_{2}] \right)^{2}/ 
(2 E[X_{2}^{2}]) $, it follows that  
$$ C^{2} \leq \frac{8}{ E[X_{2}^{2}]}B ( \theta_{1}) . $$
Combining the above inequalities with and 
$ 2 \vert D \vert \leq \vert 2D-AC \vert + \vert A \vert \vert C \vert $, 
and using (\ref{eq:bda}) and Lemmas 3.2, 3.3, 3.4, 3.5 and 4.1 
gives the following lemma.

\begin{lem} 
Assume that $(a),$ $(b),(c)$ and $(A1)$ hold. Then 
there exist $ c_{10}>0, \lambda_{*} \in (1/2,1)$ and 
$ s_{*} \in (0 ,1)$ such that, for $ \lambda \in ( \lambda_{*} ,1)$ and 
$ 0 < \vert \theta_{1} \vert \leq s_{*}$, 
$$ \vert a_{ \lambda } ( \theta_{1}) - \tilde{a}_{ \lambda }( \theta_{1}) 
\vert \leq c_{10}  \left( \vert \theta_{1} \vert ^{ \hat{\delta}} +
k_{ \lambda } ^{ \delta /4 }( \theta_{1}) \vert 
\log k_{ \lambda }( \theta_{1}) \vert \right) k_{ \lambda }
^{ -1/4 }( \theta_{1}),$$
$$ \vert b_{ \lambda } ( \theta_{1}) - \tilde{b}_{ \lambda }( \theta_{1}) 
\vert \leq c_{10}  \vert \theta_{1} \vert ^{ \tilde{\delta}} k_{ \lambda }
^{ -1/4 }( \theta_{1}), \quad \quad \quad \quad \quad
 \quad \quad \quad \quad \quad \quad $$
$$ \left\vert \frac{ b_{ \lambda }( \theta_{1}) }{ 
\sqrt{ a_{ \lambda }^{2}( \theta_{1})+ 
b_{ \lambda }^{2}( \theta_{1})}} - 
\frac{ \tilde{b}_{ \lambda }( \theta_{1}) }{ 
\sqrt{ \tilde{a}_{ \lambda }^{2}( \theta_{1})+ 
\tilde{b}_{ \lambda }^{2}( \theta_{1})}} \right\vert \leq c_{10} 
\vert \theta_{1} \vert ^{ \tilde{ \delta }} \vert \log k_{ \lambda }
( \theta_{1}) \vert ,$$
where $ \tilde{\delta}= \min \{ \hat{\delta }, \alpha \delta /2 \} $. 
Moreover, there exist $ c_{11}>0, \lambda^{*} \in (1/2,1),
s^{*} \in (0,1)$ such that, for $ \lambda \in ( \lambda^{*} ,1)$ and 
$ 0< \vert \theta_{1} \vert \leq s^{*}$, 
\begin{equation}
 \left\vert \frac{ b_{ \lambda }( \theta_{1}) }{ 
\sqrt{ a_{ \lambda }^{2}( \theta_{1})+ 
b_{ \lambda }^{2}( \theta_{1})}} - \frac{ \lambda c_{2} \vert \theta_{1}
 \vert ^{ \alpha } ( \theta_{1} / 
\vert \theta_{1} \vert ) }{ \sqrt{2} 
Q_{ \lambda }( \vert \theta_{1} \vert ^{ \alpha } )} \right\vert 
\leq c_{11}   \vert \theta_{1} \vert ^{ \delta_{0} }  ,
\label{eq:blt}
\end{equation}
where $ \delta_{0}=  \min \{ \varepsilon , \tilde{ \delta }/2 \} $.
\end{lem}

\begin{lem} 
Assume that $(a),$ $(b),(c)$ and $(A1)$ hold. Then 
there exist $ c_{12}>0, \lambda_{0} \in (1/2,1) $ and 
$ s_{0} \in (0,1)$ such that, for $ \lambda \in ( \lambda_{0} ,1)$ and 
$ 0< \vert \theta_{1} \vert \leq s_{0}$, 
\begin{equation}
 \left\vert \arcsin \frac{ b_{ \lambda }( \theta_{1}) }{ 
\sqrt{ a_{ \lambda }^{2}( \theta_{1})+ 
b_{ \lambda }^{2}( \theta_{1})}}  - 
\arcsin \frac{ \lambda c_{2} \vert \theta_{1}
 \vert ^{ \alpha } ( \theta_{1} / 
\vert \theta_{1} \vert ) }{ \sqrt{2} 
Q_{ \lambda }( \vert \theta_{1} \vert ^{ \alpha } )} \right\vert 
\leq c_{12} \vert \theta_{1} \vert ^{ \delta_{0} }
\label{eq:asas}
\end{equation}
\end{lem}
$Poof.$ By applying the Mean Value Theorem to the inverse of the sin  
function, $ \arcsin t $, 
on the closed interval with end points 
$ \lambda c_{2} \vert \theta_{1}
 \vert ^{ \alpha } ( \theta_{1} / 
\vert \theta_{1} \vert )/ (\sqrt{2} 
Q_{ \lambda }( \vert \theta_{1} \vert ^{ \alpha } )) $ 
and $ b_{ \lambda }( \theta_{1}) / 
\sqrt{ a_{ \lambda }^{2}( \theta_{1})+ 
b_{ \lambda }^{2}( \theta_{1})} $, the left-hand side of (\ref{eq:asas}) 
is equal to 
\begin{equation}
 \frac{ 1}{ \sqrt{1 -t^{2}_{0}}} \left\vert 
\frac{ b_{ \lambda }( \theta_{1}) }{ 
\sqrt{ a_{ \lambda }^{2}( \theta_{1})+ 
b_{ \lambda }^{2}( \theta_{1})}} - \frac{ \lambda c_{2} \vert \theta_{1}
 \vert ^{ \alpha } ( \theta_{1} / 
\vert \theta_{1} \vert ) }{ \sqrt{2} 
Q_{ \lambda }( \vert \theta_{1} \vert ^{ \alpha } ) } \right\vert , 
\label{eq:1bl}
\end{equation}
where $t_{0}$ is a number between $ \lambda c_{2} \vert \theta_{1}
 \vert ^{ \alpha } ( \theta_{1} / 
\vert \theta_{1} \vert )/ (\sqrt{2} 
Q_{ \lambda }( \vert \theta_{1} \vert ^{ \alpha } )) $ 
and $ b_{ \lambda }( \theta_{1}) / 
\sqrt{ a_{ \lambda }^{2}( \theta_{1})+ 
b_{ \lambda }^{2}( \theta_{1})} $.

Note that for $ \theta_{1} \ne 0$, 
\begin{equation}
\left\vert \frac{ \lambda c_{2} \vert \theta_{1}
 \vert ^{ \alpha } ( \theta_{1} / 
\vert \theta_{1} \vert ) }{ \sqrt{2} 
Q_{ \lambda }( \vert \theta_{1} \vert ^{ \alpha } ) } 
\right\vert \leq \frac{1}{ \sqrt{2}} .
\label{eq:ql}
\end{equation}
If $ \theta_{1} $ is small enough, then (\ref{eq:blt}) implies that 
$$ \left\vert \frac{ b_{ \lambda }( \theta_{1}) }{ 
\sqrt{ a_{ \lambda }^{2}( \theta_{1})+ 
b_{ \lambda }^{2}( \theta_{1})}} \right\vert \leq \frac{1}{ \sqrt{2}} + c_{11} \vert \theta_{1} 
\vert ^{ \delta_{0} } < \frac{1}{2} \left( 1+ \frac{1}{ \sqrt{2}} \right), $$
and we obtain 
$$
 \vert t_{0} \vert \leq \max \left\{ 
\left\vert \frac{ \lambda c_{2} \vert \theta_{1}
 \vert ^{ \alpha } ( \theta_{1} / 
\vert \theta_{1} \vert ) }{ \sqrt{2} 
Q_{ \lambda }( \vert \theta_{1} \vert ^{ \alpha } ) } \right\vert , 
\left\vert \frac{ b_{ \lambda }( \theta_{1}) }{ 
\sqrt{ a_{ \lambda }^{2}( \theta_{1})+ 
b_{ \lambda }^{2}( \theta_{1})}} \right\vert \right\} 
 \leq \frac{1}{2} \left( 1+ \frac{1}{ \sqrt{2}} \right) < 1.
$$
(\ref{eq:1bl}), together with (\ref{eq:blt}) and the above inequality, gives 
the desired estimate (\ref{eq:asas}).

To estimate $ \lim_{ L \rightarrow \infty } \Re [C_{L}( \lambda )] $, 
we decompose $ \Re [C_{L}( \lambda )] $. With $s_{0} >0$ being 
the same constant as in Lemma 4.3,   
$$ \Re [C_{L}( \lambda )] = 
\frac{1}{  \pi } \int_{-s_{0}}^{s_{0}} \frac{1}{ \theta_{1}} \left(
\arcsin \frac{ \lambda c_{2} \vert \theta_{1}
 \vert ^{ \alpha } ( \theta_{1} / 
\vert \theta_{1} \vert ) }{ \sqrt{2} 
Q_{ \lambda }( \vert \theta_{1} \vert ^{ \alpha } ) } \right) \ d 
\theta_{1} \quad \quad \quad \quad $$ 
$$  \quad \quad \quad \quad \quad
+ I_{1}(s_{0}, \lambda )
 + I_{2}( \lambda ,L) + I_{3}( \lambda ,L) + I_{4}( s_{0}, \lambda ) , $$
where 
$$ I_{1}(s_{0}, \lambda ) = \frac{1}{ 2\pi } \int_{-s_{0}}^{s_{0}} 
\frac{ \sin \theta_{1}}{ 1- \cos \theta_{1}} 
\left(  \arcsin  \frac{ b_{ \lambda }( \theta_{1})}
{ \sqrt{a^{2}_{ \lambda }( \theta_{1}) + b^{2}_{ \lambda }( \theta_{1})} } 
\right)   \quad \quad \quad \quad \quad \quad \quad \quad \quad \quad $$
$$ - \frac{2}{ \theta_{1}} \left(
\arcsin \frac{ \lambda c_{2} \vert \theta_{1}
 \vert ^{ \alpha } ( \theta_{1} / 
\vert \theta_{1} \vert ) }{ \sqrt{2} Q_{ \lambda }( \vert \theta_{1} \vert ^{ \alpha } )} \right) \ d \theta_{1}, $$
$$ I_{2}( \lambda ,L)= \frac{1}{ 2\pi } \int_{ - \pi }^{ \pi } 
\left(  \arcsin  \frac{ b_{ \lambda }( \theta_{1})}
{ \sqrt{a^{2}_{ \lambda }( \theta_{1}) + b^{2}_{ \lambda }( \theta_{1})} } 
\right) \sin (L \theta_{1}) \ d \theta_{1} ,
 \quad \quad \quad \quad \quad \quad \quad \quad \quad $$
$$ I_{3}( \lambda ,L) = - \frac{1}{ 2\pi } \int_{ - \pi }^{ \pi } 
\left( \arcsin  \frac{ b_{ \lambda }( \theta_{1})}
{ \sqrt{a^{2}_{ \lambda }( \theta_{1}) + b^{2}_{ \lambda }( \theta_{1})} } 
\right) \frac{ \sin \theta_{1}}{ 1- \cos \theta_{1}} 
\cos (L \theta_{1}) \ d \theta_{1}, \quad \quad \quad $$
$$  I_{4}( s_{0}, \lambda ) = \frac{1}{ 2 \pi } \int_{s_{0} < \vert \theta_{1} \vert \leq \pi }
\frac{ \sin \theta_{1}}{ 1- \cos \theta_{1}} 
\left(  \arcsin  \frac{ b_{ \lambda }( \theta_{1})}
{ \sqrt{a^{2}_{ \lambda }( \theta_{1}) + b^{2}_{ \lambda }( \theta_{1})} } 
\right) \ d \theta_{1} .
 \quad \quad \quad \quad \quad $$

Since $ \arcsin t \in [ - \pi /2 ,\pi/2 ] $, 
$ \arcsin (b_{ \lambda }( \theta_{1}) / \sqrt{ a^{2}_{ \lambda }
( \theta_{1}) + b^{2}_{ \lambda }( \theta_{1})} ) $ is 
integrable on $[- \pi , \pi]$. 
The Riemann-Lebesgue Lemma implies that 
$$ \lim_{ L \rightarrow \infty } I_{2}( \lambda ,L)=0 .$$

It is easy to verify that 
$$ \vert I_{4}( s_{0}, \lambda ) \vert \leq \frac{1}{ 2 \pi }
 \int_{s_{0} < 
\vert \theta_{1} \vert \leq \pi } \left\vert 
\frac{ \sin \theta_{1}}{ 1- \cos \theta_{1}} 
\left(  \arcsin  \frac{ b_{ \lambda }( \theta_{1})}
{ \sqrt{a^{2}_{ \lambda }( \theta_{1}) + b^{2}_{ \lambda }( \theta_{1})} } 
\right) \right\vert  
\ d \theta_{1} $$
$$ \leq  \frac{ \pi }{2 } \frac{1}{ 1 - \cos s_{0}}  . \quad 
\quad \quad \quad \quad \quad \quad \quad \quad \quad 
\quad \quad  \quad \quad $$

By applying Mean Value Theorem to 
the function $ \arcsin t $ 
on the closed interval with end points 
$ \lambda c_{2} \vert \theta_{1}
 \vert ^{ \alpha } ( \theta_{1} / 
\vert \theta_{1} \vert )/ (\sqrt{2} Q_{ \lambda }
( \vert \theta_{1} \vert ^{ \alpha } )) $ 
and $ 0$, and using (\ref{eq:ql}), 
$$ \left\vert \arcsin \frac{ \lambda c_{2} \vert \theta_{1}
 \vert ^{ \alpha } ( \theta_{1} / 
\vert \theta_{1} \vert ) }{ \sqrt{2} 
Q_{ \lambda }( \vert \theta_{1} \vert ^{ \alpha } ) } \right\vert  
 \leq \sqrt{2}  \left\vert \frac{ \lambda c_{2} \vert \theta_{1}
 \vert ^{ \alpha } ( \theta_{1} / 
\vert \theta_{1} \vert ) }{ \sqrt{2} 
Q_{ \lambda }( \vert \theta_{1} \vert ^{ \alpha } ) }
\right\vert 
 \leq \frac{ \lambda c_{2}\vert \theta_{1} \vert ^{ \alpha }}{ 1 - \lambda } .$$ 
This inequality 
and (\ref{eq:asas}) give 
\[ 
 \left\vert \left( \arcsin  \frac{ b_{ \lambda }( \theta_{1})}
{ \sqrt{a^{2}_{ \lambda }( \theta_{1}) + b^{2}_{ \lambda }( \theta_{1})} } 
\right) \frac{ \sin \theta_{1}}{ 1- \cos \theta_{1}} \right\vert 
\leq \left\{ 
\begin{array}{l}
\displaystyle{ \left( \frac{ \lambda c_{2} \vert \theta_{1} \vert ^{ \alpha} }{  1 - \lambda  } + 
c_{12} \vert \theta_{1} \vert ^{ \delta_{0}} \right) \frac{ \pi }{ \vert \theta_{1} \vert } } \quad ( 0 < \vert \theta_{1} \vert \leq s_{0} ) \\
\displaystyle{ \frac{ \pi }{2} \frac{1}{ 1- \cos s_{0}}} \quad ( s_{0} < \vert \theta_{1} \vert \leq \pi ) 
\end{array} \right.
\]
For fixed $ \lambda \in (\lambda_{0}, 1)$, the above inequality implies that 
the left-hand side of the above is integrable on $[ -\pi ,\pi]$. By the 
Riemann-Lebesgue Lemma,  
$$ \lim_{ L \rightarrow \infty } I_{3}( \lambda ,L)=0 .$$

From  (\ref{eq:asas}) and $  \arcsin t \in [ - \pi /2, \pi /2]$, the absolute 
value of the integrand in $I_{1}(s_{0}, \lambda )$ is bounded by 
$ c_{12} \pi \vert \theta_{1} \vert ^{ \delta_{0}-1} + ( \pi /2) \pi \vert \theta_{1} \vert $. Hence $I_{1}(s_{0}, \lambda )$ is 
bounded by a constant 
independent of $ \lambda \in ( \lambda_{0} ,1)$. 
We have the following estimate.

There exists a constant $c_{13}$ such that, for 
$ \lambda \in ( \lambda_{0} ,1)$,
\begin{equation}
\left\vert \lim_{ L \rightarrow \infty } \Re[ C_{L}( \lambda ) ] - 
\frac{1}{ \pi } \int_{ -s_{0}}^{s_{0}} \frac{1}{ \theta_{1}} 
\left( \arcsin \frac{ \lambda c_{2} \vert \theta_{1}
 \vert ^{ \alpha } ( \theta_{1} / 
\vert \theta_{1} \vert ) }{ \sqrt{2} 
Q_{ \lambda }( \vert \theta_{1} \vert ^{ \alpha } ) } \right) \ d 
\theta_{1} \right\vert \leq c_{13}. 
\label{eq:rcl}
\end{equation}

\begin{lem}
Assume that $(a)$, $(b)$, $(c)$ and $(A1)$ hold. Then 
there exists a constant $c_{14}>0$ such that, for 
$ \lambda \in ( \lambda_{0} ,1)$,
\begin{equation}
 \left\vert \lim_{ L \rightarrow \infty } \Re[ C_{L}( \lambda ) ] + 
\frac{2}{ \alpha \pi } \left( \arcsin  
\frac{c_{2}}{ \sqrt{2} \sqrt{ c_{1}^{2}+c_{2}^{2}+ c_{1}
\sqrt{c_{1}^{2}+c_{2}^{2}}}} \right) \log (1 - \lambda )  \right\vert 
\leq c_{14}.
\label{eq:170}
\end{equation}
\end{lem}

$Proof.$ (\ref{eq:rcl}) reduces our problem to the estimate for 
$$ I_{0}(s_{0}, \lambda ):= \frac{1}{ \pi } \int_{ -s_{0}}^{s_{0}} \frac{1}{ \theta_{1}} 
\left( \arcsin \frac{ \lambda c_{2} \vert \theta_{1}
 \vert ^{ \alpha } ( \theta_{1} / 
\vert \theta_{1} \vert ) }{ \sqrt{2} 
Q_{ \lambda }( \vert \theta_{1} \vert ^{ \alpha } ) } \right) \ d 
\theta_{1} .$$
To estimate this integral, we use the power series representation for 
$ \arcsin t $: 
$$ \arcsin t= \sum_{n=0}^{ \infty } 
\frac{ 1 \cdot 3 \cdot 5 \cdots (2n-1) }{2 \cdot 4 \cdot 6 \cdots 2n} 
\frac{t^{2n+1}}{ 2n+1} \quad ( -1 <t < 1) .$$ 
By the dominated convergence theorem, 
\begin{equation}
  I_{0}(s_{0}, \lambda ) =  \sum_{n=0}^{ \infty } 
\frac{ 1 \cdot 3 \cdot 5 \cdots (2n-1) }{2 \cdot 4 \cdot 6 \cdots 2n} 
\frac{1}{ 2n+1} \frac{1}{ \pi } 
\int_{ -s_{0}}^{s_{0}} \frac{1}{ \theta_{1}} 
\left(  
\frac{ \lambda c_{2} \vert \theta_{1}
 \vert ^{ \alpha } ( \theta_{1} / 
\vert \theta_{1} \vert ) }{ \sqrt{2} 
Q_{ \lambda }( \vert \theta_{1} \vert ^{ \alpha } ) } 
\right) ^{2n+1} \ d \theta_{1} .
\label{eq:160}
\end{equation}
We notice that the function 
$$ \frac{1}{ \theta_{1}} \left( \frac{  \lambda c_{2} 
\vert \theta_{1} \vert ^{ \alpha } 
( \theta_{1} / \vert \theta_{1} \vert ) }{ \sqrt{2} Q_{ \lambda }
( \vert \theta_{1} \vert ^{ \alpha })} \right) ^{2n+1} $$ 
that appears on the right-hand side of the above equality 
is an even function, and 
perform the change of variables 
$ s =  \theta_{1}  ^{ \alpha } $, so that 
$$ \int_{ -s_{0}}^{s_{0}} \frac{1}{ \theta_{1}} 
\left(  
\frac{ \lambda c_{2} \vert \theta_{1}
 \vert ^{ \alpha } ( \theta_{1} / 
\vert \theta_{1} \vert ) }{ \sqrt{2} 
Q_{ \lambda }( \vert \theta_{1} \vert ^{ \alpha } ) } 
\right) ^{2n+1} \ d \theta_{1} = 2 \int_{0}^{s_{0}} \frac{1}{ \theta_{1}} 
\left(  
\frac{ \lambda c_{2} \vert \theta_{1}
 \vert ^{ \alpha } ( \theta_{1} / 
\vert \theta_{1} \vert ) }{ \sqrt{2} 
Q_{ \lambda }( \vert \theta_{1} \vert ^{ \alpha } ) } 
\right) ^{2n+1} \ d \theta_{1} \quad \quad \quad  $$ 
\begin{equation}
 \quad \quad \quad \quad \quad \quad \quad \quad \quad \quad \quad 
\quad  
= \frac{2}{ \alpha } \left( \frac{ \lambda c_{2}}{ \sqrt{2}} \right) ^{2n+1} 
 \int_{0}^{ s_{0}^{ \alpha } } 
\frac{ s^{2n}}{ ( Q_{ \lambda }( s))^{2n+1}} \ ds . 
\label{eq:150}
\end{equation}

Put 
$$ Q_{ \lambda }^{(1)}(s) = \left\{ \lambda^{2} (c_{1}^{2}+c_{2}^{2} + c_{1} 
\sqrt{c_{1}^{2}+c_{2}^{2}} )s^{2}+ (1 - \lambda )\lambda (3c_{1}+ \sqrt{ 
c_{1}^{2}+c_{2}^{2}} )s + 2( 1 - \lambda )^{2} \right\} ^{1/2}  $$
and 
$$ Q_{ \lambda }^{(2)}(s) = \left\{ \lambda^{2} (c_{1}^{2}+c_{2}^{2} + c_{1} 
\sqrt{c_{1}^{2}+c_{2}^{2}} )s^{2} \right. 
\quad \quad \quad \quad \quad \quad \quad \quad \quad \quad \quad 
 \quad \quad \quad \quad \quad \quad \quad \quad $$
$$  \quad \quad \quad \quad
+ (1 -\lambda ) \lambda \left( 2c_{1}+ 
\sqrt{c_{1}^{2}+c_{2}^{2}}+ \frac{c_{1}^{2}}{ \sqrt{c_{1}^{2}+c_{2}^{2}}} 
\right) s
+ \left. (1- \lambda )^{2} \left( 1+ \frac{c_{1}}
{ \sqrt{c_{1}^{2}+c_{2}^{2}}} \right)
\right\} ^{1/2} .$$
$ Q_{ \lambda }^{(1)}(s)$ is obtained by replacing $ \sqrt{(1- \lambda + \lambda c_{1}s)^{2}+ ( \lambda c_{2}s)^{2}} $ in the definition of $Q_{ \lambda }(s) $ with $ \lambda \sqrt{c_{1}^{2}+c_{2}^{2}} s + 1- \lambda $. 
$ Q_{ \lambda }^{(2)}(s)$ is obtained by replacing $ \sqrt{(1- \lambda + \lambda c_{1}s)^{2}+ ( \lambda c_{2}s)^{2}} $ in the definition of $Q_{ \lambda }(s) $ with $ \lambda \sqrt{c_{1}^{2}+c_{2}^{2}} s + (1- \lambda)(c_{1}/ 
\sqrt{c_{1}^{2}+c_{2}^{2}}) $.

From these inequalities, 
$$ 
\lambda \sqrt{ c_{1}^{2}+c_{2}^{2}} s + (1 - \lambda ) 
\frac{c_{1}}{ \sqrt{ c_{1}^{2}+c_{2}^{2}}} \leq 
\sqrt{(1- \lambda + \lambda c_{1}s)^{2}+ ( \lambda c_{2}s)^{2}} 
\leq   \lambda \sqrt{c_{1}^{2}+c_{2}^{2}} s + 1- \lambda $$ 
for $ s>0$, so we obtain 
\begin{equation}
 \int_{0}^{s_{0}^{ \alpha }} \frac{s^{2n}}{
(Q^{(1)}_{ \lambda }(s)) ^{2n+1}} \ ds \leq 
\int_{0}^{s_{0}^{ \alpha }} \frac{s^{2n}}{
(Q_{ \lambda }(s)) ^{2n+1}} \ ds \leq 
\int_{0}^{s_{0}^{ \alpha }} \frac{s^{2n}}{
(Q_{ \lambda }^{(2)}(s)) ^{2n+1}} \ ds .
\label{eq:140}
\end{equation}

For $ a,b,c \in (0, \infty) $ with $ 4ac -b^{2}>0$, we next evaluate 
the integral 
$$ I_{n}^{ \alpha }(s_{0}) := \int_{0}^{ s_{0}^{ \alpha }} \frac{s^{2n}}{ 
(as^{2}+bs+c)^{(2n+1)/2} } \ ds \quad ( n \in \{0 \} \cup {\Bbb N}).$$
In the case when $n=0$, 
$$ I_{0}^{ \alpha }(s_{0})= \frac{1}{ \sqrt{a}} \left[ \log (2as+ b + 2 
\sqrt{a(as^{2}+bs+c)} \right]_{0}^{s_{0}^{ \alpha } } . $$
By setting $a= \lambda ^{2}(c_{1}^{2}+c_{2}^{2}+ 
c_{1} \sqrt{c_{1}^{2}+c_{2}^{2}} ) , b= (1- \lambda ) \lambda (3c_{1}+ 
\sqrt{c_{1}^{2}+c_{2}^{2}}),$ $ c=2(1 - \lambda )^{2}$, this identity 
yields the following estimate.

There exists a constant $c_{15} >0$ such that, for $ \lambda \in (1/2,1)$, 
\begin{equation}
 \left\vert \frac{ \lambda c_{2}}{ \sqrt{2}}
 \int_{0}^{s_{0}^{ \alpha} } \frac{1}
{ Q_{ \lambda }^{(1)}(s)}  \ ds + \frac{c_{2}}{ 
\sqrt{2} \sqrt{c_{1}^{2}+c_{2}^{2}+ 
c_{1} \sqrt{c_{1}^{2}+c_{2}^{2}}}} \log (1 - \lambda ) \right\vert \leq 
c_{15} . 
\label{eq:ioes}
\end{equation}

In the case when $ n \in {\Bbb N}$, we make the substitutions 
$$ u= \frac{1}{(2n-1)a} \frac{1}{ (as^{2}+bs+ c)^{(2n-1)/2}}, \quad \quad 
dv = (2n-1)s^{2n-2}, $$ 
$$ du= - \frac{1}{2a} \frac{2as+b}{ (as^{2}+bs+ c)^{(2n+1)/2}}, \quad 
\quad \quad v= s^{2n-1} . \quad \quad \quad $$
The formula for integration by parts gives 
\begin{equation}
 I_{n}^{ \alpha }(s_{0}) = \frac{1}{ a} I_{n-1}^{ \alpha }(s_{0}) + e_{n}, 
\label{eq:inita}
\end{equation}
where 
$$ e_{n}= - \frac{1}{(2n-1)a} \frac{s_{0}^{ (2n-1) \alpha }}
{ (as_{0}^{2 \alpha } +b s_{0}^{ \alpha }+c)^{(2n-1)/2}} 
- \frac{b}{2a} \int_{0}^{s_{0}^{ \alpha }} \frac{s^{2n-1}}
{(as^{2}+bs+c)^{(2n+1)/2}} \ ds . $$

Since $a,b,c$ are positive, satisfying the inequality $ 4ac-b^{2} >0$,
$$ \vert e_{n} \vert \leq \frac{1}{ a^{(2n+1)/2}} + \frac{b}{ 2a^{(2n+1)/2}}
 \int_{0}^{ s_{0}^{ \alpha }} \frac{1}{ as^{2}+bs+c} \ ds $$ 
\begin{equation}
 \leq 
 \frac{1}{ a^{(2n+1)/2}} \left( 1+ \frac{ \pi }{2} 
\frac{b}{ \sqrt{4ac - b^{2}}} \right) . \quad \quad \quad 
\label{eq:enera}
\end{equation}

By using (\ref{eq:inita}) repeatedly, we obtain 
$$
 I_{n}^{ \alpha }(s_{0}) = \frac{1}{a^{n}} 
I_{0}^{ \alpha }(s_{0}) + \sum_{j=0}^{n-1} \frac{1}{a^{j}}e_{n-j} .
$$

With the help of the above relation and (\ref{eq:enera}), 
and by setting $a= \lambda ^{2}(c_{1}^{2}+c_{2}^{2}+ 
c_{1} \sqrt{c_{1}^{2}+c_{2}^{2}} ) , b= (1- \lambda ) \lambda (3c_{1}+ 
\sqrt{c_{1}^{2}+c_{2}^{2}}),$ $ c=2(1 - \lambda )^{2}$, (\ref{eq:ioes}) 
implies the following estimate for $k=1$.

There exists a constant $c_{16}>0$ such that, for $ \lambda \in (1/2,1)$ and 
$ n \in \{ 0 \} \cup {\Bbb N} $, 
$$ \left\vert \left( \frac{ \lambda c_{2}}{ \sqrt{2}} \right)^{ 2n+1} 
\int_{0}^{ s_{0}^{ \alpha } } \frac{ s^{2n}}{ (Q_{ \lambda }
^{(k)}(s) )^{2n+1} } \ ds + \left(
\frac{c_{2}}{ \sqrt{2} \sqrt{c_{1}^{2}+c_{2}^{2}+ 
c_{1} \sqrt{c_{1}^{2}+c_{2}^{2}}}} \right)^{2n+1}
\log (1 - \lambda ) \right\vert $$
\begin{equation}
 \leq c_{16}( n+1) \left( \frac{1}{ \sqrt{2}}
\right)^{2n} ,  \quad (k=1,2) .  \quad \quad
 \quad \quad \quad \quad \quad \quad \quad \quad \quad \quad 
 \quad \quad \quad \quad  
\label{eq:qnnn}
\end{equation}
The same argument as in the proof of the estimate for $k=1$ 
in (\ref{eq:qnnn}) gives the estimate for $k=2$ in (\ref{eq:qnnn}).

By (\ref{eq:140}), (\ref{eq:qnnn}) also holds when the integral of 
$ s^{2n}/(Q^{(k)}_{ \lambda }(s) )^{2n+1}$ over $[0,s^{ \alpha }_{0} ]$ 
with respect to $ s $ is replaced by the integral of 
$ s^{2n}/(Q_{ \lambda }(s) )^{2n+1}$ over $[0,s^{ \alpha }_{0} ]$ 
with respect to $ s $. This estimate implies (\ref{eq:170}) in view of 
(\ref{eq:150}) and (\ref{eq:160}).

We recall (\ref{eq:180}) and (\ref{eq:190}). It immediate from 
Lemma 4.4 that, for $ \lambda \in ( \lambda_{0} , 1 ) ,$
$$  (1- \lambda )^{-2 \beta }
e^{-c_{14}} \leq \frac{ P_{0} \{ T_{ \lambda }  < \tau_{V_{-}} \} }
{ P_{0} \{ T_{ \lambda }  < \tau_{V_{+}} \}  }  \leq (1- \lambda )^{-2 \beta }
e^{c_{14}}, $$ 
where $ \beta $ is the same constant as in Theorem 1.1. By combining 
the above inequalities, Lemma 2.1 and Lemma 2.2, we have 
the following estimate. 
\begin{pro}
Assume that $(a)$, $(b)$, $(c)$ and $(A1)$ hold. Then 
there exist constants $c_{17}>0,c_{18}>0$ and $ \lambda^{*}_{0} \in (1/2 ,1)$ 
such that, for $ \lambda \in ( \lambda^{*}_{0} ,1)$,
$$ c_{17} (1- \lambda )^{- \beta+(1/4) }
 \leq  P_{0} \{ T_{ \lambda }  < \tau_{V_{-}} \} = 
(1- \lambda ) \sum_{n=0}^{ \infty }
\lambda^{n} P_{0} \{ \tau_{V_{-}} >n \} 
 \leq c_{18} (1- \lambda )^{- \beta+(1/4) } .$$
\end{pro}
This asymptotic behavior of $ P_{0} \{ T_{ \lambda }  < \tau_{V_{-}} \} $ 
implies (\ref{eq:200}) by the Tauberian theorem (see [5, Theorem 2.4.3]). 
The proof of Theorem 1.1 is complete.

The proofs of Corollaries 1.2 and 1.3 are similar to that of Theorem 1.1 and 
are thus omitted.

\section{Appendix}

We calculate 
$$ \Re \left[ \frac{1}{ 2 \pi} \int_{ - r_{0} /2 }^{ r_{0}/2 } \frac{1}{ 
\varphi _{ \lambda } (  \theta_{1}, \theta_{2} )  } \ d \theta_{2}
\right]  = \frac{1}{ 2 \pi} \int_{ - r_{0} /2 }^{ r_{0}/2 } \Re \left[
 \frac{1}{ 
\varphi _{ \lambda } (  \theta_{1}, \theta_{2} )  } \right] \ d \theta_{2} , $$
where $r_{0}$ is the constant in Lemma 3.1.

Recall that 
$$ A= \frac{ E[ ( \sin \theta_{1}X_{1}) X_{2}] }{ (1/2) E[X_{2}^{2} ] } ,
\quad B= \frac{ 1 - \lambda E[ \cos \theta_{1}X_{1}] }
{( \lambda /2) E[X_{2}^{2} ] } ,$$
$$ C= \frac{ E[ ( -1 + \cos  \theta_{1}X_{1} )X_{2} ]}{ (1/2) E[X_{2}^{2} ] } ,
\quad D = \frac{ E[ \sin \theta_{1}X_{1} ] }{  (1/2) E[X_{2}^{2} ] }, $$
$$ K= ( 4B -A^{2}+C^{2})^{2}+ 4 (2D-AC)^{2} $$
and 
$$  \varphi _{ \lambda } 
(  \theta_{1}, \theta_{2} )  =
 \frac{ \lambda }{2} E[X_{2}^{2} ] \{ \theta_{2}^{2} + A \theta _{2} +B
 - i (C \theta_{2}+D ) \}  . $$

Put 
$$ a_{ \pm }= A \pm \sqrt{ H } , $$
\[ b_{ \pm } = \left\{
\begin{array}{ll}
 \displaystyle{ \frac{1}{4} \left\{ A^{2} +C^{2} + \sqrt{K} \pm 
  \frac{ 2AH - 2C(2D-AC)}
{ \sqrt{ H }} \right\} }  & 
2D - AC \ne 0 \\
\displaystyle{ \frac{1}{2} \{ 2B +C^{2} \pm 
\sqrt{4B-A^{2}+C^{2}} \vert C \vert \} }  & 
2D - AC = 0,
\end{array} \right.  \]
where
$$ H = \frac{1}{2} (-4B + A^{2}- C^{2}  + \sqrt{K} ) .$$
Then $a_{ \pm },b_{ \pm }$ are real numbers, since $4B-A^{2}-C^{2}>0$, 
and 
$$ ( \theta_{2}^{2} + A \theta_{2} +B)^{2} + ( C \theta_{2} + 
D)^{2} = ( \theta_{2}^{2} +a_{+} \theta_{2} + b_{+} )
( \theta_{2}^{2} +a_{-} \theta_{2} + b_{-}) .$$

In the case $ 2D-AC \ne 0$, 
it is easy to verify that 
\begin{equation}
 (a_{+}-a_{-})(a_{+}b_{-}-a_{-}b_{+}) + 
(b_{-} - b_{+})^{2} = \frac{1}{2} ( -4B +A^{2}+C^{2} + \sqrt{K}) 
\sqrt{K} \ne 0 
\label{eq:abb}
\end{equation}
using 
\begin{equation}
(b_{-} -b_{+})^{2} = \left\{ 
-   \frac{ AH - C(2D-AC)}
{ \sqrt{ H }} \right\}^{2}
 = \frac{1}{4} (A^{2} + C^{2} + \sqrt{K} )^{2} - 
4(B^{2}+D^{2}) . \quad \quad \quad 
\label{eq:3}
\end{equation}
The partial-fraction decomposition of 
$ \Re [1/ \varphi _{ \lambda } ( \theta_{1}, \theta_{2}) ] $ is 
$$ \Re \left[ \frac{1}{ \varphi _{ \lambda } ( \theta_{1}, \theta_{2})} \right]
 = \frac{ \Re[ \varphi _{ \lambda } ( \theta_{1}, \theta_{2}) ] }{
\vert \varphi _{ \lambda } ( \theta_{1}, \theta_{2}) \vert ^{2}} 
= \frac{1}{ ( \lambda /2)E[X_{2}^{2}] }
 \left(
\frac{  F \theta_{2}+G }
{ \theta_{2}^{2} + a_{+} \theta_{2} + b_{+}} + 
\frac{  -F \theta_{2}+ I}
{ \theta_{2}^{2} + a_{-} \theta_{2} + b_{-}} \right), $$
where 
$$ F=  \frac{ -(a_{+}b_{-}-a_{-}b_{+}) +A(b_{-}-b_{+}) + 
B(a_{+}-a_{-})}{(a_{+}-a_{-})(a_{+}b_{-}-a_{-}b_{+}) + 
(b_{-}-b_{+})^{2}} , \quad \quad \quad \quad \quad \quad
 \quad $$
$$ G= \frac{ -b_{+}( b_{-}-b_{+}) -A(a_{+}-a_{-})b_{+} +B
a_{+}(a_{+}-a_{-})+B(b_{-}-b_{+}) }{(a_{+}-a_{-})(a_{+}b_{-}-a_{-}b_{+}) + 
(b_{-}-b_{+})^{2}}, $$
$$ I= \frac{ b_{-}(b_{-}-b_{+}) +A (a_{+}-a_{-})b_{-}-B 
a_{-}(a_{+}-a_{-})-B(b_{-}-b_{+})}{(a_{+}-a_{-})(a_{+}b_{-}-a_{-}b_{+}) + 
(b_{-}-b_{+})^{2}}.$$
From the definition of $a_{ \pm} $ and $ b_{ \pm } $ 
$$ 4b_{ \pm } - a_{ \pm }^{2} = J_{1} \pm J_{2} ,$$
where 
$$ J_{1}= \frac{1}{2} (4B -A^{2} +3C^{2} + \sqrt{K} ),\quad
J_{2} = - \frac{ 2C(2D-AC)}{ \sqrt{H}} .$$

From (\ref{eq:3}),  
$ 4B -A^{2} -C^{2} > 0$ implies 
\begin{equation}
J_{1}^{2}-J_{2}^{2}= \frac{1}{4} ( 4B -A^{2}-C^{2} + \sqrt{K} )^{2} >0.
\label{eq:j}
\end{equation}
Combining the above inequality with $ J_{1}>0$, we have 
$$ \vert J_{2} \vert < J_{1} $$ 
and  
\begin{equation}
4b_{ \pm } - a_{ \pm }^{2} >0 . 
\label{eq:ba}
\end{equation}

We turn now to the integration. In the case $2D -AC \ne 0$, 
$$
\frac{ \lambda E[X_{2}^{2}]}{2} 
 \int_{ - r_{0}/2}^{ r_{0}/2 } \Re \left[ 
 \frac{ 1}{ \varphi_{ \lambda } ( \theta_{1} , \theta_{2} )}  \right] 
d \theta_{2}  \quad \quad \quad \quad \quad \quad \quad
 \quad \quad \quad \quad \quad \quad \quad \quad \quad \quad \quad
 \quad \quad \quad $$
$$ = \int_{ - r_{0}/2 }^{ r_{0}/2 } \frac{ F \theta_{2} +G }{ \theta_{2}^{2} + a_{+} 
\theta_{2} +b_{+}} + \frac{ - F \theta_{2} + I }
{ \theta_{2}^{2} + a_{-} 
\theta_{2} +b_{-}} \ d \theta_{2}  \quad \quad \quad \quad \quad 
 \quad \quad \quad  \quad \quad \quad$$
$$ = \frac{ F}{2} 
\left[ \log \frac{ \theta_{2} ^{2} +a_{+} \theta_{2} +b_{+} }{ \theta_{2} ^{2} + a_{-} \theta_{2} + b_{-}} \right] _{ - r_{0}/2 }^{ r_{0}/2} + \frac{ -a_{+}F +2G }{ \sqrt{ 4b_{+} -a_{+}^{2}}} 
\left[ \arctan \frac{ 2 \theta_{2} +a_{+}}{ \sqrt{ 4b_{+} -a_{+}^{2}}} \right]_{ -r_{0}/2 }^{ r_{0}/2 } 
$$ 
\begin{equation}
  +  \frac{ a_{-}F +2 I }{ \sqrt{ 4b_{-} -a_{-}^{2}}} 
\left[ \arctan \frac{ 2 \theta_{2} +a_{-}}{ \sqrt{ 4b_{-} -a_{-}^{2}}} \right]_{ - r_{0}/2 }^{ r_{0}/2 }. 
 \quad \quad \quad \quad \quad \quad \quad \quad \quad
\label{eq:mr}
\end{equation}

In the case $2D -AC = 0$, we obtain 
$$
\frac{ \lambda E[X_{2}^{2}]}{2} 
 \int_{ - r_{0}/2}^{ r_{0}/2 } \Re \left[ 
 \frac{ 1}{ \varphi_{ \lambda } ( \theta_{1} , \theta_{2} )} \right]
d \theta_{2}  \quad \quad \quad \quad \quad \quad \quad
 \quad \quad \quad \quad \quad \quad \quad \quad \quad \quad \quad
 \quad \quad \quad $$
$$
= \frac{1}{ \sqrt{4B -A^{2} +C^{2}} } \left[ \arctan \frac{ 2 \theta_{2} +A}
{ \sqrt{ 4B -A^{2} +C^{2}} + \vert C \vert }
\right]_{ - r_{0} /2}^{ r_{0} /2 }  \quad \quad  
$$
$$
 + \frac{1}{ \sqrt{4B -A^{2} +C^{2}} } \left[ \arctan \frac{ 2 \theta_{2} +A}
{ \sqrt{ 4B -A^{2} + C^{2}} - \vert C \vert }
\right]_{ - r_{0} /2}^{ r_{0} /2 } . 
$$

{\bf Proof of Lemma 3.2.} We will prove (\ref{eq:l2}) in the case $ 2D- AC \ne 0$; the proof 
of (\ref{eq:l2}) in the case $ 2D-AC=0 $ is similar and is omitted.

Assume $ 2D- AC \ne 0$. In view of (\ref{eq:3}) and the definition of $a_{ \pm }$, $a_{ \pm }$ and $b_{ \pm } $ tend to 
$0$ as $ \theta _{1} \rightarrow 0$ and $ \lambda \uparrow 1$. 
With the help of 
$$  F= - \frac{ \sqrt{H}}{ \sqrt{K}} \quad  \mbox{ and } \quad 
 \vert \log (1+x ) \vert \leq 2 \vert x \vert \quad
 (\vert x \vert \leq 1/2), $$
the first term on the right-hand side of 
the last equality in (\ref{eq:mr}) is bounded by a constant.

A simple calculation gives
$$ - a_{+} F +2 G = \frac{1}{2 \sqrt{K}} \left\{ 4B -A^{2}+C^{2} + \sqrt{K} - 
\frac{ 2C(2D-AC) }{ \sqrt{H}} \right\} ,$$
$$ a_{-} F +2 I = \frac{1}{2 \sqrt{K}} \left\{ 4B -A^{2}+C^{2} + \sqrt{K} + \frac{2C(2D-AC)}{ \sqrt{H}} \right\} .$$
We recall that 
$$ J_{2}= - \frac{2C(2D-AC)}{ \sqrt{H}} . $$
Since  
$$ \vert J_{2} \vert < J_{1} = 
\frac{1}{2} ( 4B -A^{2} +3C^{2}  + \sqrt{K} ),$$
as shown in the proof of (\ref{eq:ba}), 
$ \vert - a_{+} F +2 G \vert $ is bounded by a constant. 
The difference between the second term on the right-hand side of 
the last equality in (\ref{eq:mr}) 
and 
$$ \frac{ -a_{+}F +2G }{ \sqrt{ 4b_{+} - a_{+}^{2}}} \pi 
=  \frac{ \pi }{2 \sqrt{ J_{1} + J_{2}} \sqrt{K}} \left\{ 4B -A^{2}+C^{2} + 
\sqrt{K} - \frac{ 2C(2D-AC) }{ \sqrt{H}} \right\}$$
is bounded by a constant. Here we used the inequality 
$$ 0 \leq \frac{ \pi }{2} - \arctan x \leq \frac{1}{x} \quad ( x >0).$$

Similarly, the difference between the third term on the 
right-hand side of the last equality in (\ref{eq:mr}) 
and 
$$ \frac{ a_{-}F +2 I }{ \sqrt{ 4b_{-} - a_{-}^{2}}} \pi 
=  \frac{ \pi }{2 \sqrt{ J_{1} - J_{2}} \sqrt{K}} \left\{ 4B -A^{2}+C^{2} + 
\sqrt{K} + \frac{ 2C(2D-AC) }{ \sqrt{H}} \right\}$$
is bounded by a constant. Put  
$$ \bar{a}_{ \lambda }( \theta_{1}) =
\frac{1}{ 2 \lambda E[X_{2}^{2}] \sqrt{K} } \left( 
\frac{ J + J_{2}}{ \sqrt{ J_{1}+J_{2}} } 
+ \frac{ J - J_{2}}{  \sqrt{ J_{1}-J_{2}} } \right),
$$
where 
$$ J= 4B-A^{2}+C^{2} + \sqrt{K}.$$
Then the difference between $ (1/2 \pi) \int_{ - r_{0}/2 }^{ r_{0}/2 } \Re [ 
1/ \varphi_{ \lambda } ( \theta_{1} ,\theta_{2})] d \theta_{2} $ and $ \bar{a}_{ \lambda }( \theta_{1}) $ is bounded by a constant.

It remains to show that 
\begin{equation}
 \bar{a}_{ \lambda }( \theta_{1}) = 
\tilde{a}_{ \lambda }( \theta_{1}) := \frac{ \sqrt{2}}{ \lambda 
E[X_{2}^{2}] \sqrt{K}} \sqrt{ 4B-A^{2}+C^{2}+ \sqrt{K}}. 
\label{eq:a}
\end{equation} 
Since $ 4B -A^{2}-C^{2} >0 $, 
$$ \vert J_{2} \vert < J_{1} \leq J, $$ 
which implies that $ \bar{a}_{ \lambda }( \theta_{1}) \geq 0$.
(The first inequality above has already been shown in the proof of 
(\ref{eq:ba}).)
We will compute 
$ \left( \bar{a}_{ \lambda }( \theta_{1}) \right)^{2} $. By (\ref{eq:j}) and 
$ J= 2J_{1}- 2C^{2} $, 
$$ \left\{ (J+J_{2}) \sqrt{ J_{1}-J_{2}} + (J-J_{2}) \sqrt{ J_{1}+J_{2}}
\right\}^{2} 
 \quad \quad \quad \quad \quad \quad \quad \quad \quad
 \quad \quad \quad \quad \quad \quad \quad $$ 
$$ = 2 \left( J_{1}+ \sqrt{ J_{1}^{2}-J_{2}^{2}} 
\right)J^{2}+ 2 \left( J_{1}-2J- 
\sqrt{ J_{1}^{2}-J_{2}^{2}} \right)J_{2}^{2} 
 \quad \quad \quad \quad \quad \quad \quad \quad \quad 
 \quad $$
$$ = 8(4B -A^{2}+C^{2}+ \sqrt{K})(J_{1}-C^{2})^{2} 
- 4(4B-A^{2}+ \sqrt{K})J_{2}^{2} 
\quad \quad \quad \quad \quad \quad \quad \ $$
$$
= 8(4B -A^{2}+C^{2}+ \sqrt{K})(J_{1}^{2}-J_{2}^{2}) + 
8(4B -A^{2}+C^{2}+ \sqrt{K})(-2J_{1}+C^{2})C^{2} 
$$
\begin{equation}
+ 4(4B-A^{2}+2C^{2}+ \sqrt{K})J_{2}^{2} . \quad 
 \quad \quad \quad \quad \quad \quad \quad \quad
 \quad \quad \quad \quad \quad \quad \quad \quad
\label{eq:jj}
\end{equation}
In view of 
$$ J_{2}^{2} =J_{1}^{2} - \frac{1}{4}(4B-A^{2}-C^{2}+ \sqrt{K})^{2} 
 = 2(4B- A^{2}+C^{2}+ \sqrt{K})C^{2}, $$ 
which follows from (\ref{eq:j}), the sum of the second and third terms 
on the right-hand side of the last equality in (\ref{eq:jj}) vanish. 
Hence 
$$ \left( \bar{a}_{ \lambda }( \theta_{1}) \right)^{2} 
= \left( \tilde{a}_{ \lambda }( \theta_{1}) \right)^{2} .$$ 
Since $
\bar{a}_{ \lambda }( \theta_{1})$ and 
$ \tilde{a}_{ \lambda }( \theta_{1})$ are nonnegative functions, 
the above relation yields (\ref{eq:a}).

To show Lemma 3.3, we require the following lemmas. 
 
\begin{lem}
Assume $(a)$, $(b)$ and $(c)$ hold. Then, for $ \lambda \in (1/2,1) $ and 
$ \theta_{1} \in [ - \pi ,\pi ] $ with $ \vert \theta _{1} \vert \leq 
r_{0}/2 $, 
$$ (i) \quad  \frac{1}{ 2 \pi } \int_{ - r_{0}/2 }^{ r_{0}/2 } 
\frac{ \theta_{2}^{2} }{ \vert \varphi _{ \lambda } (  \theta_{1}, \theta_{2} )  \vert ^{2} } \ 
d \theta_{2} \leq   \frac{ 16} {  c_{*} E[X_{2}^{2}]}  
K^{ -1/4} , 
 \quad \quad \quad \quad \quad \quad \quad \quad \quad \quad \quad
 \quad \quad \quad $$ 
$$ (ii) \quad \frac{1}{ 2 \pi } \int_{ - r_{0}/2 }^{ r_{0}/2 } 
\frac{ 1 }{ \vert \varphi _{ \lambda } (  \theta_{1}, \theta_{2} ) 
 \vert ^{2} } \ d \theta_{2} \leq  
\frac{1}{ 2 \pi } \int_{ - \infty }^{ \infty }
\frac{ 1 }{ \vert \varphi _{ \lambda } (  \theta_{1}, \theta_{2} ) 
 \vert ^{2} } \ d \theta_{2} \leq  
\frac{ 64}{ ( E[X_{2}^{2}] )^{2}} K ^{-3/4}. \quad $$
\end{lem}

$Proof.$ By (\ref{eq:22}), we obtain that, for 
$ \vert \theta_{1} \vert \leq r_{0}/2 $,  
$$ \int_{- r_{0}/2}^{r_{0}/2} \frac{ \theta_{2}^{2}}{ \vert \varphi_{ \lambda }
( \theta_{1}, \theta_{2}) \vert ^{2} } \ d \theta_{2} \leq 
\frac{4}{c_{*}} \int_{-r_{0}/2}^{r_{0}/2} \frac{ \Re [ 
\varphi_{ \lambda }( \theta_{1} , \theta_{2} )  ] }{ \vert \varphi_{ \lambda }
( \theta_{1}, \theta_{2}) \vert ^{2} } \ d \theta_{2} \leq \frac{4}{c_{*}} \int_{ - \infty }^{ \infty} \Re \left[ \frac{1}{\varphi_{\lambda }( \theta_{1} , \theta_{2} ) } \right] \ 
d \theta_{2} . $$ Assume that $ 2D-AC \ne 0$. 
In view of (\ref{eq:mr}), the definition of $ \tilde{a}_{ \lambda }
 ( \theta_{1})$ and $4B -A^{2}+C^{2} \leq \sqrt{K}$ imply that 
$$ \int_{ - \infty }^{ \infty } \Re \left[ \frac{1}{ \varphi_{\lambda }
( \theta_{1} , \theta_{2} ) }\right] \ d \theta_{2} 
= 2 \pi \tilde{a}_{ \lambda }
 ( \theta_{1}) \leq \frac{ 4 \pi} { \lambda E[X_{2}^{2}]}  
K^{ -1/4} .  $$
The above two inequalities give the desired estimate (i) in the case 
$2D-AC \ne 0$. The proof of (i) in the case $ 2D-AC =0 $ is similar and 
is omitted.

Since $ 1/ \vert \varphi_{ \lambda }( \theta_{1}, \theta_{2}) 
\vert ^{2} \geq 0$, the first inequality of (ii) is clear.
To obtain the last inequality of (ii), we calculate the integral of $ 1/ \vert \varphi_{ \lambda }( \theta_{1}, \theta_{2}) \vert ^{2}$ over $( - \infty , \infty )$ with respect to 
$ \theta_{2}$. Assume $ 2D-AC \ne 0$. 
The partial-fraction decomposition of $ 1/ \vert \varphi_{ \lambda }( \theta_{1}, \theta_{2}) \vert ^{2} $ is 
$$ \frac{1}{ \vert \varphi_{ \lambda }( \theta_{1}, \theta_{2}) \vert ^{2}} = 
\frac{4}{ ( \lambda E[X_{2}^{2}])^{2}} \left( \frac{ \bar{F} \theta_{2} 
+ \bar{G} }{ \theta_{2}^{2} + a_{+} \theta_{2} + b_{+} } +
\frac{ - \bar{F} \theta_{2} 
+ \bar{I} }{ \theta_{2}^{2} + a_{-} \theta_{2} + b_{-} } 
\right), $$
where 
$$ \bar{F}= \frac{a_{+}-a_{-} }{ (a_{+}-a_{-})(a_{+}b_{-} -a_{-}b_{+} ) + 
(b_{-}-b_{+})^{2} } , $$
$$ \bar{G}= \frac{ a_{+} (a_{+} -a_{-} ) + (b_{-}-b_{+} )}
{ (a_{+}-a_{-})(a_{+}b_{-} -a_{-}b_{+} ) + 
(b_{-}-b_{+})^{2} } , $$
$$ \bar{I}= \frac{ -a_{-}(a_{+} -a_{-}) -(b_{-}-b_{+}) }
{ (a_{+}-a_{-})(a_{+}b_{-} -a_{-}b_{+} ) + 
(b_{-}-b_{+})^{2} } . $$
This gives the integral 
$$ \frac{( \lambda E[X_{2}^{2}])^{2}}{4} \int_{- \infty}^{ \infty } 
\frac{1}{ \vert \varphi_{ \lambda }( \theta_{1}, \theta_{2}) \vert ^{2} } \ d 
\theta_{2} = \left( \frac{ -a_{+} \bar{F} + 2 \bar{G}}{ \sqrt{ 4b_{+} -a_{+}^{2}}} + \frac{ a_{-} \bar{F} + 2 \bar{I}}{ \sqrt{ 4b_{-} -a_{-}^{2}}} 
\right) \pi.$$
A simple calculation and (\ref{eq:abb}) show that the right-hand side of 
the above equality is equal to 
$$ c_{\lambda }( \theta_{1} ):= \left\{ 
\frac{2H -J_{2} }{ \sqrt{J_{1}+J_{2}}} 
+ \frac{ 2H +J_{2} }{ \sqrt{J_{1}-J_{2}}} \right\} 
 \times \frac{ 2 \pi}{ (-4B+A^{2}+C^{2} + \sqrt{K}) \sqrt{K}}.  $$
We compute $(c_{\lambda }( \theta_{1} ))^{2}$ by using (\ref{eq:j}). 
From this 
computation and $ c_{\lambda }( \theta_{1} ) \geq 0$, we have 

$$ c_{\lambda }( \theta_{1} )= \frac{ 4 \sqrt{2} \pi }{ \sqrt{K}}
\frac{ \sqrt{4B-A^{2}+C^{2}+ \sqrt{K}}}{4B-A^{2}-C^{2} + \sqrt{K}}.
$$
This implies the desired estimate (ii) in the case 
$2D-AC \ne 0$. The proof of (ii) in the case $ 2D-AC =0 $ is similar and 
is omitted.

\begin{lem}
Assume that $(a)$, $(b)$ and $(c)$ hold. Then 
there exists a constant $c>0$ such that, 
for $ \lambda \in (1/2,1) $ and 
$ \theta_{1} \in [ - \pi ,\pi ] $ with $ \vert \theta _{1} \vert \leq 
r_{0}/2 $,  
$$ (iii) \quad \frac{1}{2 \pi } \int_{-r_{0}/2}^{r_{0}/2} 
\frac{ \vert \theta_{2} \vert ^{2+ \delta } }
{ \vert \varphi _{ \lambda } (  \theta_{1}, \theta_{2} )  \vert ^{2} } \ 
d \theta_{2} \leq c K^{ (-1+ \delta )/4 } ( 1 + \vert \log K \vert ), 
\quad \quad \quad \quad \quad \quad  \quad \quad \quad \quad \quad \quad $$
$$ (iv) \quad  \frac{1}{2 \pi } \int_{-r_{0}/2}^{r_{0}/2} 
\frac{ \vert \theta_{2} \vert ^{ \gamma } }
{ \vert \varphi _{ \lambda } (  \theta_{1}, \theta_{2} )  \vert ^{2} } \ 
d \theta_{2} \leq c K^{ (-3+ \gamma )/4 } \quad 
( \gamma = \delta ,1, 1+ \delta ).
 \quad \quad \quad \quad \quad \quad \quad \quad \quad $$

\end{lem}

$Proof.$ We prove (iii). Since $ \delta \in (0,1)$, 
\begin{equation}
 \vert \theta_{2} \vert ^{2+ \delta } \leq \left\{ 
\begin{array}{ll}
K^{(2+ \delta )/4} & ( \vert \theta_{2} \vert \leq K^{1/4} ) \\
 \vert \theta_{2} \vert ^{3} K^{(-1+ \delta )/4} & 
(\vert \theta_{2} \vert > K^{1/4} ).
\end{array} 
\right.  \label{eq:t2}
\end{equation}
We use (\ref{eq:t2}) to obtain that 
the left-hand side of (iii) is bounded by 
\begin{equation}
 \frac{K^{ (2+ \delta )/4}}{ 2 \pi} \int_{- r_{0}/2 }^{ r_{0}/2 } 
\frac{1}{ \vert \varphi_{ \lambda } ( \theta_{1}, \theta_{2} ) \vert ^{2} } 
\ d \theta_{2} 
+  \frac{K^{(-1+ \delta )/4}}{ 2 \pi } \int_{ K^{1/4} < \vert \theta_{2} 
\vert \leq r_{0}/2} 
\frac{ \vert \theta_{2} \vert ^{3}}
{ \vert \varphi_{ \lambda } ( \theta_{1}, \theta_{2} ) \vert ^{2} } 
 \ d \theta_{2} .
\label{eq:t2d}
\end{equation}
By (ii) in Lemma 5.1, the first term of (\ref{eq:t2d}) is 
bounded by a constant multiple of 
$K^{(-1+ \delta )/4}$. From (\ref{eq:22}), we replace 
$ \vert \varphi_{ \lambda } ( \theta_{1}, \theta_{2} )  \vert^{2} $ with 
$ c_{*}^{2} \vert \theta_{2} \vert ^{4} /16 $ in the second term 
of (\ref{eq:t2d}). Then the second term of (\ref{eq:t2d}) is 
bounded by a constant multiple of 
$ K^{(-1+ \delta )/4} (1+ \vert \log K \vert ).$ Thus we conclude that 
(iii) holds.

To prove (iv), 
we use the inequality 
\[
 \vert \theta_{2} \vert ^{ \gamma } \leq \left\{ 
\begin{array}{ll}
K^{ \gamma /4} & ( \vert \theta_{2} \vert \leq K^{1/4} ) \\
 \vert \theta_{2} \vert ^{2} K^{( \gamma -2 )/4} & 
(\vert \theta_{2} \vert > K^{1/4} )
\end{array} 
\right. 
\]
instead of (\ref{eq:t2}), and the result follows by (i) and (ii) in Lemma 5.1.

{\bf Proof of Lemma 3.3.} With $r_{0} >0$, 
$$ a_{\lambda }( \theta_{1}) - \frac{1}{2 \pi } \int_{-r_{0}/2}^{r_{0}/2}
\Re \left[\frac{1}{ \varphi_{ \lambda }( \theta_{1}, \theta_{2}) } \right]
\ d \theta_{2}
= \frac{1}{2 \pi } \int_{ -r_{0}/2} ^{r_{0}/2}  \chi_{ \lambda }
 ( \theta_{1}, \theta_{2} ) \ d \theta_{2} \quad 
\quad  \quad \quad \quad \quad \quad \quad \quad \quad \quad $$
\begin{equation} 
 \quad \quad \quad \quad \quad \quad \quad \quad \quad \quad \quad
 \quad
+ \frac{1}{2 \pi } \int_{ r_{0}/2 < \vert \theta_{2} \vert \leq \pi } 
\frac{ \Re[ 1- \lambda \phi ( \theta_{1} ,\theta_{2} ) ]}{ \vert 
1- \lambda \phi ( \theta_{1} ,\theta_{2} ) \vert ^{2} } \ d \theta_{2},  
\label{eq:r1}
\end{equation}
where
$$ \chi_{ \lambda } ( \theta_{1}, \theta_{2} ) 
 = \frac{\Re[ 1- \lambda \phi ( \theta_{1} ,\theta_{2} ) ]}{ \vert 
1- \lambda \phi ( \theta_{1} ,\theta_{2} ) \vert ^{2} } - 
\frac{ \Re [ \varphi_{ \lambda }( \theta_{1}, \theta_{2} )]}{
\vert \varphi_{ \lambda }( \theta_{1}, \theta_{2} ) \vert ^{2}}.  $$
(\ref{eq:2}) implies that, for $ \lambda \in (1/2,1)$ and 
$( \theta_{1}, \theta_{2} ) \in [ - \pi ,\pi ] 
\times [ - \pi ,\pi ] $ 
with $ \vert \theta_{2} \vert \geq r_{0}/2$, 
$$ \left\vert \frac{\Re[ 1- \lambda \phi ( \theta_{1} ,\theta_{2} ) ]}{ \vert 
1- \lambda \phi ( \theta_{1} ,\theta_{2} ) \vert ^{2} } \right\vert 
\leq \frac{2}{ \vert 
1- \lambda \phi ( \theta_{1} ,\theta_{2} ) \vert ^{2}}   
\leq \frac{128}{c_{*}^{2}r_{0}^{4}} .$$  
Hence the second term on the right-hand side of (\ref{eq:r1}) is bounded by a 
constant.

To estimate the first term on the right-hand side of (\ref{eq:r1}), 
we consider $ \chi_{ \lambda } ( \theta_{1}, \theta_{2} ).$ 
By applying the estimates (\ref{eq:1}) and 
$$ \frac{1}{ \vert 1 - \lambda \phi ( \theta_{1} ,\theta_{2} ) \vert } 
\leq \frac{ c_{3}( \pi^{ \hat{ \delta }} + \pi^{ \delta } )+1 }{  
\vert \varphi_{ \lambda } ( \theta_{1} ,\theta_{2} ) \vert } \quad  
( \lambda \in (1/2,1), 
(\theta_{1} ,\theta_{2}) \in [- \pi,\pi ] \times [- \pi,\pi ]), $$
which follows from (\ref{eq:vp}), 
we deduce that 
$$ \vert  \chi_{ \lambda } ( \theta_{1}, \theta_{2} ) \vert 
\leq \vert \Re[ 1 - \lambda \phi ( \theta_{1} ,\theta_{2} ) ] 
- \Re [ \varphi_{ \lambda }( \theta_{1} ,\theta_{2} ) ] \vert 
\frac{1}{ \vert 1 - \lambda \phi ( \theta_{1} ,\theta_{2} ) \vert ^{2} } 
\quad \quad \quad \quad \quad \quad \quad \quad \quad \quad 
 \quad \quad$$
$$  \quad \quad \quad \quad \quad
+ \vert \Re [ \varphi_{ \lambda }( \theta_{1} ,\theta_{2} ) ] \vert
\frac{ \vert 1 - \lambda \phi ( \theta_{1} ,\theta_{2} ) 
- \varphi_{ \lambda }( \theta_{1} ,\theta_{2} ) \vert }
{ \vert 1 - \lambda \phi ( \theta_{1} ,\theta_{2} ) \vert 
\vert \varphi_{ \lambda }( \theta_{1} ,\theta_{2} ) \vert } 
\left( \frac{1}{ \vert 1 - \lambda \phi ( \theta_{1} ,\theta_{2} ) \vert } + 
\frac{1}{ \vert \varphi_{ \lambda }( \theta_{1} ,\theta_{2} ) \vert } \right)
$$
\begin{equation}
 \leq \hat{c} \frac{( \vert \theta_{1} \vert ^{ \hat{ \delta }} + \vert \theta_{2} \vert ^{\delta }) \theta_{2}^{2} }
{ \vert \varphi_{ \lambda }( \theta_{1} ,\theta_{2} ) \vert^{2}}
, \quad \quad \quad  \quad \quad \quad \quad \quad \quad \quad \quad
 \quad \quad \quad \quad \quad \quad \quad \quad \quad \quad
\label{eq:cii}
\end{equation}
where $ \hat{c} $ is a constant independent of $ \lambda \in (1/2,1)$ 
and $ ( \theta_{1} ,\theta_{2} ) \in [ - \pi ,\pi ] \times [ - \pi ,\pi ].$

The proof completed by showing the following estimates. 
The integral of 
$ \vert \theta_{2} \vert ^{2}/ \vert 
\varphi_{ \lambda }( \theta_{1} ,\theta_{2} ) \vert ^{2} $ 
appearing on the 
right-hand side of the above inequality over $[ -r_{0}/2 ,r_{0}/2]$ 
with respect to $ \theta_{2} $ is bounded by a constant multiple of 
$K^{-1/4} .$ Moreover, the integral of 
$ \vert \theta_{2} \vert ^{2+ \delta }/ \vert 
\varphi_{ \lambda }( \theta_{1} ,\theta_{2} ) \vert ^{2} $ 
appearing on the 
right-hand side of the above inequality over $[ -r_{0}/2 ,r_{0}/2]$ 
with respect to $ \theta_{2} $ is bounded by a constant multiple of 
$ K^{ (-1+ \delta )/4 } \vert \log K \vert $. 
These estimates come from (i) in Lemma 5.1 and (iii) in Lemma 5.2. 
Notice that $K$ tends to zero as $ \theta_{1} \rightarrow 0$ and 
$ r \uparrow 1 $, which implies that $ 1 \leq \vert \log K \vert $.

{\bf Proof of Lemma 3.4.} (\ref{eq:33}) can be proved by the same procedure as 
in Lemma 3.2. 
We only give a sketch of the proof of (\ref{eq:33}) in 
the case $2D-AC \ne 0.$

Assume that $ 2D-AC \ne 0$. We calculate 
$$ \Im \left[ \frac{1}{2 \pi} \int_{-r_{0}/2}^{r_{0}/2} 
\frac{1}{ \varphi _{ \lambda } ( \theta_{1} ,\theta_{2} )} \ d \theta_{2} 
\right] = \frac{1}{2 \pi} \int_{-r_{0}/2}^{r_{0}/2} 
\Im \left[ \frac{1}{ \varphi _{ \lambda } ( \theta_{1} ,\theta_{2} )} 
 \right] \ d \theta_{2} . $$
The partial-fraction decomposition of 
$ \Im [ 1/ \varphi _{ \lambda } ( \theta_{1} ,\theta_{2} )] $ 
is 
$$ \Im \left[ \frac{1}{ \varphi _{ \lambda } ( \theta_{1} ,\theta_{2} )} 
 \right] = \frac{ - \Im [ \varphi _{ \lambda } ( \theta_{1} ,\theta_{2} ) ]}
{ \vert \varphi _{ \lambda } ( \theta_{1} ,\theta_{2} ) \vert^{2}} 
= \frac{1}{ ( \lambda /2 ) E[X_{2}^{2}] } \left( 
 \frac{ \tilde{F} \theta_{2} + \tilde{G} }{ \theta_{2}^{2} + a_{+} 
\theta_{2} + b_{+} } + \frac{ - \tilde{F} \theta_{2} + \tilde{I} }{ \theta_{2}^{2} + a_{-} 
\theta_{2} + b_{-} } \right), $$
where 
$$ \tilde{F}= \frac{ C(b_{-}-b_{+} )+ D(a_{+}-a_{-} )}
{ (a_{+}-a_{-})(a_{+}b_{-}- a_{-}b_{+} ) + (b_{-} -b_{+} )^{2}} ,
\quad \quad \quad \quad \quad $$
$$ \tilde{G}= \frac{ -C(a_{+} -a_{-})b_{+} + D a_{+}(a_{+}-a_{-} ) 
+ D(b_{-}-b_{+} )}{ (a_{+}-a_{-})(a_{+}b_{-}- a_{-}b_{+} ) + 
(b_{-} -b_{+} )^{2}} ,$$
$$ \tilde{I}= \frac{ C(a_{+} -a_{-})b_{-} - D a_{-}(a_{+}-a_{-} ) 
- D(b_{-}-b_{+} )}{ (a_{+}-a_{-})(a_{+}b_{-}- a_{-}b_{+} ) + 
(b_{-} -b_{+} )^{2}} .$$
This gives the integral 
$$ \frac{ \lambda E[X_{2}^{2}]}{2} \int_{- r_{0}/2}^{ r_{0}/2 } 
\Im \left[ \frac{1}{  \varphi_{ \lambda }( \theta_{1}, \theta_{2})  } 
\right] \ d \theta_{2} \quad \quad \quad \quad \quad \quad \quad
 \quad \quad \quad \quad \quad \quad \quad \quad \quad \quad \quad
 \quad \quad \quad \quad  $$
$$ 
= \frac{ \tilde{F} }{2} \left[ \log \frac{ \theta_{2}^{2} +a_{+} \theta_{2} 
+b_{+}}{\theta_{2}^{2} +a_{-} \theta_{2} +b_{-} } \right]_{ -r_{0}/2}^{r_{0}/2} + \frac{ -a_{+} \tilde{F}+ 2 \tilde{G}}{ \sqrt{ 4b_{+}-a_{+}^{2}}} 
\left[ \arctan \frac{ 2 \theta_{2}+ a_{+}}{ \sqrt{ 4b_{+} -a_{+}^{2}}} 
\right]_{ - r_{0}/2}^{r_{0}/2} 
$$
\begin{equation}
+ \frac{ a_{-} \tilde{F}+ 2 \tilde{I}}{ \sqrt{ 4b_{-}-a_{-}^{2}}} 
\left[ \arctan \frac{ 2 \theta_{2}+ a_{-}}{ \sqrt{ 4b_{-} -a_{-}^{2}}} 
\right]_{ - r_{0}/2}^{r_{0}/2} . \quad \quad \quad \quad \quad \quad 
\quad \quad \quad 
\label{eq:limb}
\end{equation}
Using 
$$ \tilde{F} = \frac{ 2D-AC }{ \sqrt{K} \sqrt{H}}, \quad \quad 
\vert \log (1+x) \vert \leq 2 \vert x \vert \quad 
( \vert x \vert \leq 1/2 ), $$
and 
\begin{equation}
 (2D-AC)^{2}=  \frac{1}{2}(4B -A^{2} +C^{2} + \sqrt{K})H , 
\label{eq:j2}
\end{equation}
the first term on the right-hand side of (\ref{eq:limb}) is bounded by 
a constant multiple of $ \vert 2D- AC \vert / \sqrt{K} + \vert C \vert .$ 
A simple calculation gives 
$$ -a_{+} \tilde{F}+ 2 \tilde{G}= \frac{ 2D- AC }{ \sqrt{K}} - \frac{ C \sqrt{H}}{ \sqrt{K}}, \quad \quad 
  a_{-} \tilde{F}+ 2 \tilde{I}= \frac{ 2D- AC }{ \sqrt{K}} + \frac{ C \sqrt{H}}{ \sqrt{K}} $$
by $(\ref{eq:j2})$ above.
The difference between the second term on the right-hand side of 
(\ref{eq:limb}) and 
$$  \frac{ -a_{+} \tilde{F}+ 2 \tilde{G} }{ \sqrt{4b_{+} -a_{+}^{2}}} \pi 
= \frac{ \pi }{ \sqrt{ J_{1} +J_{2}} \sqrt{K}} (2D- AC - C \sqrt{H} )$$
is bounded by a constant multiple of $ \vert 2D - AC \vert / \sqrt{K} + \vert C \vert /K^{1/4}.$ Similarly, the difference between the third term on the 
right-hand side of (\ref{eq:limb}) and 
$$  \frac{ a_{-} \tilde{F}+ 2 \tilde{I} }{ \sqrt{4b_{-} -a_{-}^{2}}} \pi 
= \frac{ \pi }{ \sqrt{ J_{1} -J_{2}} \sqrt{K}} (2D- AC + C \sqrt{H} )$$
is bounded by a constant multiple of $ \vert 2D - AC \vert / \sqrt{K} + \vert C \vert /K^{1/4}.$ 

Put 
$$ \bar{b}_{ \lambda }( \theta_{1}) 
= \frac{1}{ \lambda E[X_{2}^{2}] \sqrt{K} } \left( 
\frac{ 2D-AC -C \sqrt{H}}{ \sqrt{ J_{1}+J_{2}} } 
+ \frac{ 2D-AC + C \sqrt{H}}{  \sqrt{ J_{1}-J_{2}} } \right) .
$$
Then the difference between $ (1/2 \pi) \int_{ - r_{0}/2 }^{ r_{0}/2 } \Im [ 1/ \varphi_{ \lambda } ( \theta_{1} ,\theta_{2})] d \theta_{2} $ and 
$ \bar{b}_{ \lambda }( \theta_{1}) $ is 
bounded by a constant multiple of $ \vert 2D - AC \vert / \sqrt{K} + \vert C \vert /K^{1/4}.$

It remains to show that 
$$ \bar{b}_{ \lambda } ( \theta_{1} )= \tilde{b}_{ \lambda }( \theta_{1})
:= \frac{ \sqrt{2}}{ \lambda E[X_{2}^{2}] \sqrt{K}} 
\frac{2(2D-AC)}{ \sqrt{ 4B-A^{2}+C^{2}+ \sqrt{K}}}.$$
By using (\ref{eq:j}) and noticing 
$(\ref{eq:j2})$, we have 
\begin{equation}
\bar{b}^{2}_{ \lambda }( \theta_{1} )=
\frac{ 4H}{ ( \lambda E[X_{2}^{2}] )^{2} K} 
= \tilde{b}^{2}_{ \lambda }( \theta_{1}) . \label{eq:bhb}
\end{equation}
From the inequality 
$$ 0 \leq C^{2} H \leq \frac{1}{2} (4B-A^{2}+C^{2}+ \sqrt{K} )H 
= (2D-AC)^{2} ,$$ 
the sign of $ \tilde{b}_{ \lambda }( \theta_{1} ) $ is the same as 
that of $ \bar{b}_{ \lambda }( \theta_{1} ) $. 
(\ref{eq:bhb}) yields 
$ \bar{b}_{ \lambda } ( \theta_{1} )= \tilde{b}_{ \lambda }( \theta_{1}).$

{\bf Proof of Lemma 3.5.} Let $ \hat{b}_{ \lambda }( \theta_{1} ) $ be the imaginary part of 
$$ \frac{1}{2 \pi } \int_{ - r_{0}/2 }^{ r_{0}/2 } \frac{1}{ \varphi
_{ \lambda } ( \theta_{1} ,\theta_{2}) } \ d \theta_{2} . $$
By performing the change of variables $ \theta_{2}= - \tilde{\theta}_{2}$ 
and using the symmetry of trigonometric functions, 
$ b_{ \lambda }( \theta_{1}) =-b_{ \lambda }( -
\theta_{1}) $ and $ \hat{b}_{ \lambda }( \theta_{1} ) = 
- \hat{b}_{ \lambda }( - \theta_{1} )  $. In view of this, 
$$ b_{ \lambda }( \theta_{1}) - \hat{b}_{ \lambda }( \theta_{1} ) 
= \frac{1}{2} \{ b_{ \lambda }( \theta_{1}) - b_{ \lambda }(- \theta_{1}) 
- ( \hat{b}_{ \lambda }( \theta_{1} )  - \hat{b}_{ \lambda }( - \theta_{1} ) )
\} \quad \quad \quad \quad \quad \quad $$
$$ \quad \quad  \quad \quad 
= \frac{1}{ 4 \pi } \int_{-r_{0}/2}^{r_{0}/2} 
 \chi_{ \lambda }^{(1)} ( \theta_{1}, \theta_{2}) \ d \theta_{2} + 
\frac{1}{ 4 \pi } \int_{-r_{0}/2}^{r_{0}/2}
 \chi_{ \lambda }^{(2)} ( \theta_{1}, \theta_{2}) \ d \theta_{2} $$
\begin{equation}
 \quad \quad  \quad \quad +  
\frac{1}{ 4 \pi } \int_{ r_{0}/2 \leq \vert \theta_{2} \vert \leq \pi }
\chi_{ \lambda }^{(3)} ( \theta_{1}, \theta_{2}) \ d \theta_{2}, 
\quad \quad \quad \quad \quad \quad \quad 
\label{eq:b-b}
\end{equation}
where 
$$ \chi_{ \lambda }^{(1)} ( \theta_{1}, \theta_{2}) 
 = \lambda E[ ( \sin \theta_{1}X_{1} ) 
( \cos \theta_{2}X_{2} ) ] \left( \frac{ 1}{ \vert 1 - \lambda \phi 
(  \theta_{1} ,\theta_{2} ) \vert^{2} }  + \frac{ 1}{ \vert 1 - \lambda \phi 
( - \theta_{1} ,\theta_{2} ) \vert^{2} } \right) $$
$$ - \lambda E[  \sin \theta_{1}X_{1} ] \left( \frac{1}{ \vert \varphi
_{ \lambda} ( \theta_{1} ,\theta_{2} ) \vert ^{2}} + 
\frac{1}{ \vert \varphi_{ \lambda} 
 ( - \theta_{1} ,\theta_{2} ) \vert ^{2}} \right), $$
$$  \chi_{ \lambda }^{(2)} ( \theta_{1}, \theta_{2}) =  
\lambda E[ ( \cos \theta_{1}X_{1})( \sin \theta_{2} X_{2})] 
\left( \frac{ 1}{ \vert 1 - \lambda \phi 
(  \theta_{1} ,\theta_{2} ) \vert^{2} }  - \frac{ 1}{ \vert 1 - \lambda \phi 
( - \theta_{1} ,\theta_{2} ) \vert^{2} } \right) $$
$$ \quad \quad \quad
 - \lambda E[(-1+ \cos \theta_{1} X_{1}) X_{2} ] \theta_{2} 
\left( \frac{1}{ \vert \varphi
_{ \lambda} ( \theta_{1} ,\theta_{2} ) \vert ^{2}} - 
\frac{1}{ \vert \varphi_{ \lambda} 
 ( - \theta_{1} ,\theta_{2} ) \vert ^{2}} \right)
$$
and 
$$ \chi_{ \lambda }^{(3)} ( \theta_{1}, \theta_{2}) = 
\frac{ \lambda E[ \sin ( \theta_{1}X_{1}+ \theta_{2}X_{2})]}
{ \vert 1 - \lambda \phi 
(  \theta_{1} ,\theta_{2} ) \vert ^{2}} - 
\frac{ \lambda E[ \sin ( - \theta_{1}X_{1}+ \theta_{2}X_{2})]}
{ \vert 1 - \lambda \phi 
( - \theta_{1} ,\theta_{2} ) \vert ^{2}}. 
\quad \quad \quad \quad \quad $$

The inequality $ E[ \vert \sin \theta_{1} X_{1} \vert ] \leq E[ \vert X_{1} \vert ^{ \delta} ] \vert \theta_{1} \vert ^{ \delta } $ and (\ref{eq:2}) 
imply that, for 
$ \lambda \in (1/2,1)$ and $( \theta_{1}, \theta_{2} ) \in [ - \pi ,\pi ] 
\times [ - \pi ,\pi ] $ 
with $ \vert \theta_{2} \vert \geq r_{0}/2 $, 
$$ \vert \chi_{ \lambda } ^{(3)} ( \theta_{1}, \theta_{2})  \vert 
\leq \frac{\vert E[ ( \sin \theta_{1}X_{1})( \cos \theta_{2}X_{2})] \vert}
{ \vert 1 - \lambda \phi 
(  \theta_{1} ,\theta_{2} ) \vert ^{2} } + 
\frac{\vert E[ ( \sin \theta_{1}X_{1})( \cos \theta_{2}X_{2})] \vert}
{\vert 1 - \lambda \phi 
( - \theta_{1} ,\theta_{2} ) \vert ^{2}} \quad 
\quad \quad \quad \quad \quad \quad \quad $$
$$  \quad \quad \quad \quad \quad
+ \vert E[ ( \cos \theta_{1} X_{1} )( \sin \theta_{2}X_{2}) ] \vert 
\left\vert \frac{1}{ \vert 1 - \lambda \phi 
(  \theta_{1} ,\theta_{2} ) \vert ^{2} } - \frac{1}{\vert 1 - \lambda \phi 
( - \theta_{1} ,\theta_{2} ) \vert ^{2}} \right\vert $$
$$ \leq \frac{ E[\vert X_{1} \vert ^{ \delta} ] \vert \theta_{1} \vert ^{ \delta } }{\vert 1 - \lambda \phi 
(  \theta_{1} ,\theta_{2} ) \vert ^{2} } + 
\frac{ E[\vert X_{1} \vert ^{ \delta} ] \vert \theta_{1} \vert ^{ \delta } }
{\vert 1 - \lambda \phi 
( - \theta_{1} ,\theta_{2} ) \vert ^{2}} 
 \quad \quad \quad \quad \quad \quad \quad \quad $$
$$  \quad \quad \quad + \vert \phi 
(  \theta_{1} ,\theta_{2} ) - \phi 
( - \theta_{1} ,\theta_{2} ) \vert 
\frac{ \vert 1 - \lambda \phi ( - \theta_{1} ,\theta_{2} ) \vert 
+ \vert 1 - \lambda \phi (  \theta_{1} ,\theta_{2} ) \vert }
{ \vert 1 - \lambda \phi (  \theta_{1} ,\theta_{2} ) \vert ^{2}
\vert 1 - \lambda \phi ( - \theta_{1} ,\theta_{2} ) \vert^{2} } $$
$$ \leq  \frac{128 E[\vert X_{1} \vert ^{ \delta} ] }
{ c_{*}^{2}r_{0}^{4}} \left( 
1+ \frac{ 32}{c_{*}r_{0}^{2}} \right) \vert \theta_{1}
 \vert ^{ \delta }.
\quad \quad \quad \quad \quad \quad \quad \quad \quad \quad 
 \quad \quad$$
Hence the third term on the right-hand side of 
the last equality in (\ref{eq:b-b}) is bounded 
by a constant multiple of $ \vert \theta_{1} \vert ^{ \delta }$.

To estimate the first and second term on the right-hand side of 
the last equality in (\ref{eq:b-b}), 
we consider $ \chi_{ \lambda }^{(1)} ( \theta_{1}, \theta_{2}) $ and 
$ \chi_{ \lambda }^{(2)} ( \theta_{1}, \theta_{2}) $. 
With the help of 
$$ \vert E[ ( \sin \theta_{1}X_{1} )( \cos \theta_{2} X_{2}) ] - 
E[ \sin \theta_{1}X_{1} ] \vert 
\leq E[ \vert X_{1} \vert ^{ \hat{ \delta}} X_{2}^{2} ] \vert 
\theta_{1} \vert ^{ \hat{ \delta}} \theta_{2}^{2} $$ 
and (\ref{eq:2}), 
the same argument as shown in (\ref{eq:cii}) implies that 
$$ \vert \chi_{ \lambda }^{(1)} ( \theta_{1} ,\theta_{2})  \vert 
\leq  
\frac{\vert E[( \sin \theta_{1}X_{1}) ( \cos \theta_{2} X_{2})] - 
E[ \sin \theta_{1}X_{1} ] \vert}
{ \vert 1 - \lambda \phi (  \theta_{1} ,\theta_{2} ) 
\vert ^{2} } 
 \quad \quad \quad \quad \quad \quad \quad \quad 
 \quad \quad \quad \quad \quad \quad $$
$$ + \frac{\vert E[( \sin \theta_{1}X_{1}) ( \cos \theta_{2} X_{2})] - 
E[ \sin \theta_{1}X_{1} ] \vert}{ \vert 1 - \lambda \phi 
( - \theta_{1} ,\theta_{2} ) \vert^{2} } 
\quad \quad \quad \quad \quad \quad $$
$$ 
+ \vert E[ \sin \theta_{1} X_{1}] \vert  
\frac{ \vert 1 - \lambda \phi ( \theta_{1} ,\theta_{2} ) 
- \varphi_{ \lambda }( \theta_{1} ,\theta_{2} ) \vert }
{ \vert 1 - \lambda \phi ( \theta_{1} ,\theta_{2} ) \vert 
\vert \varphi_{ \lambda }( \theta_{1} ,\theta_{2} ) \vert } 
\quad \quad \quad \quad \quad $$
$$ \times 
\left( \frac{1}{ \vert 1 - \lambda \phi ( \theta_{1} ,\theta_{2} ) \vert } + 
\frac{1}{ \vert \varphi_{ \lambda }( \theta_{1} ,\theta_{2} ) \vert } \right) 
 \quad \quad \quad \quad \quad \quad $$
$$ \quad \quad \quad 
+ \vert E[ \sin \theta_{1} X_{1}] \vert  
\frac{ \vert 1 - \lambda \phi (- \theta_{1} ,\theta_{2} ) 
- \varphi_{ \lambda }( - \theta_{1} ,\theta_{2} ) \vert }
{ \vert 1 - \lambda \phi ( - \theta_{1} ,\theta_{2} ) \vert 
\vert \varphi_{ \lambda }( - \theta_{1} ,\theta_{2} ) \vert } 
\quad \quad \quad \quad \quad \quad $$
$$ \times 
\left( \frac{1}{ \vert 1 - \lambda \phi ( - \theta_{1} ,\theta_{2} ) \vert } + 
\frac{1}{ \vert \varphi_{ \lambda }( - \theta_{1} ,\theta_{2} ) \vert } \right)
 \quad \quad \quad \quad $$
\begin{equation} 
\quad \quad 
 \leq \bar{c}_{1} \{ \vert D \vert ( \vert \theta_{1} \vert ^{ \hat{ \delta}} + \vert \theta_{2} \vert ^{ \delta }) + 
\vert \theta_{1} \vert ^{ \hat{ \delta}}  \theta_{2} ^{2} \} 
\left( \frac{1}{ \vert \varphi_{ \lambda }( \theta_{1} ,\theta_{2} )
 \vert ^{2} } 
+ \frac{1}{\vert \varphi_{ \lambda }( - \theta_{1} ,\theta_{2} ) \vert ^{2} } 
\right) , \label{eq:ch1}
\end{equation}
where $ \bar{c}_{1} $ is a 
constant independent of $ \lambda \in (1/2,1)$ 
and $ ( \theta_{1} ,\theta_{2} ) \in [ - \pi ,\pi ] \times [ - \pi ,\pi ]$.

By performing the change of variables $ \theta_{2}= - \tilde{\theta}_{2}$ 
and using the symmetry of trigonometric functions,
\begin{equation}
 \int_{-r_{0}/2}^{r_{0}/2}  \frac{ \vert \theta_{2} \vert ^{ \gamma }}
{\vert \varphi_{ \lambda }(  - \theta_{1} ,\theta_{2} ) \vert ^{2} } \ 
d \theta_{2}=
\int_{-r_{0}/2}^{r_{0}/2}  \frac{ \vert \tilde{\theta}_{2} \vert ^{ \gamma }}
{\vert \varphi_{ \lambda }(   \theta_{1} , \tilde{\theta}_{2} ) \vert ^{2} } \ 
d \tilde{\theta}_{2} , \quad ( \gamma \geq 0). \label{eq:sy}
\end{equation}
Taking note of (\ref{eq:sy}), it follows from (\ref{eq:ch1}), Lemma 5.1 and (iv) in Lemma 5.2 that 
the first term on the right-hand side of 
the last equality in (\ref{eq:b-b}) is bounded by a constant multiple of 
$   \vert \theta_{1} \vert ^{ \hat{\delta}} K^{-1/4} + 
\vert D \vert ( \vert \theta_{1} \vert ^{ \hat{\delta}} + K^{ \delta /4} ) 
K^{-3/4}  .$

From conditions (b) and (c), we obtain 
$$ \vert E[ (\cos \theta_{1}X_{1}) ( \sin \theta_{2}X_{2})  ] - 
E[ (-1+ \cos \theta_{1}X_{1}) X_{2} ] \theta_{2} \vert \quad \quad \quad 
\quad \quad \quad \quad $$
$$ = \vert E[ (\cos \theta_{1}X_{1}) ( \sin \theta_{2}X_{2})  ] - 
E[ (\cos \theta_{1}X_{1}) X_{2} ] \theta_{2} \vert \leq  E [ \vert X_{2} \vert 
^{ 2+ \delta } ] \vert \theta_{2} \vert ^{2 + \delta } $$
and 
$$ \vert  \phi ( \theta_{1} ,\theta_{2} )-  \phi (- \theta_{1} ,\theta_{2} ) 
\vert \leq 2 \vert E[ \sin \theta_{1}X_{1}] \vert + 2
\vert E[ (\sin \theta_{1}X_{1} ) X_{2}] \vert \vert \theta_{2} \vert $$
$$ \quad \quad \quad \quad 
+ 2 E[ \vert X_{1} \vert ^{ \hat{\delta}}X_{2}^{2}] 
\vert \theta_{1} \vert ^{ \hat{\delta}} \theta_{2}^{2}. $$
Using the above two inequalities and (\ref{eq:2}), the same argument as given in (\ref{eq:cii}) implies that 
$$ \vert \chi_{ \lambda }^{(2)} ( \theta_{1} ,\theta_{2}) \vert \leq 
\frac{ \vert 
 E[( \cos \theta_{1}X_{1}) ( \sin \theta_{2}X_{2}) ] - 
E[ (-1+ \cos \theta_{1}X_{1}) X_{2} ] \theta_{2} \vert }
{ \vert 1 - \lambda \phi ( \theta_{1} ,\theta_{2} ) \vert 
\vert 1 - \lambda \phi ( - \theta_{1} ,\theta_{2} ) \vert} 
\quad \quad \quad \quad \quad \quad $$
$$ \quad \quad \quad \quad \quad \quad 
\times \vert 
\phi ( \theta_{1} ,\theta_{2} ) - \phi ( - \theta_{1} ,\theta_{2} ) \vert 
\left( \frac{1}{ \vert 1 - \lambda \phi ( \theta_{1} ,\theta_{2} ) \vert} + 
\frac{1}{ \vert 1 - \lambda \phi ( - \theta_{1} ,\theta_{2} ) \vert } \right)$$
$$ \quad \quad 
+ \vert E[(-1+ \cos \theta_{1}X_{1}) X_{2} ] \theta_{2} \vert 
\frac{ \vert 1 - \lambda \phi ( \theta_{1} ,\theta_{2} ) 
- \varphi_{ \lambda }( \theta_{1} ,\theta_{2} ) \vert }
{ \vert 1 - \lambda \phi ( \theta_{1} ,\theta_{2} ) \vert 
\vert \varphi_{ \lambda }( \theta_{1} ,\theta_{2} ) \vert } $$
$$ \times \left( \frac{1}{ \vert 1 - \lambda \phi ( \theta_{1} ,\theta_{2} ) 
\vert } + 
\frac{1}{ \vert \varphi_{ \lambda }( \theta_{1} ,\theta_{2} ) \vert } 
\right)  \quad \quad \quad \quad \quad $$
$$ \quad \quad \quad \quad 
+ \vert E[(-1+ \cos \theta_{1}X_{1}) X_{2} ] \theta_{2} \vert 
\frac{ \vert 1 - \lambda \phi (- \theta_{1} ,\theta_{2} ) 
- \varphi_{ \lambda }( - \theta_{1} ,\theta_{2} ) \vert }
{ \vert 1 - \lambda \phi ( - \theta_{1} ,\theta_{2} ) \vert 
\vert \varphi_{ \lambda }( - \theta_{1} ,\theta_{2} ) \vert } $$
$$ \times 
\left( \frac{1}{ \vert 1 - \lambda \phi ( - \theta_{1} ,\theta_{2} ) \vert } + 
\frac{1}{ \vert \varphi_{ \lambda }( - \theta_{1} ,\theta_{2} ) \vert } 
\right) \quad \quad $$
$$ \quad \leq  \bar{c }_{2}
 \{  \vert D \vert \vert \theta_{2} \vert ^{  \delta }+
 ( \vert A \vert   + \vert C \vert ) 
\vert \theta_{2} \vert ^{ 1+ \delta } + ( \vert C \vert \vert \theta_{2} \vert 
 + \vert \theta_{2} \vert ^{2} ) 
\vert \theta_{1} \vert ^{ \hat{ \delta}} \} $$
\begin{equation}
 \times \left(  
\frac{1}{ \vert \varphi_{ \lambda }(  \theta_{1} ,\theta_{2} ) \vert ^{2}} +
\frac{1}{ \vert \varphi_{ \lambda }( - \theta_{1} ,\theta_{2} ) \vert ^{2} }
\right) , \quad \quad \label{eq:ch2}
\end{equation}
where $ \bar{c }_{2} $ is a constant 
independent of $ \lambda \in (1/2,1)$ 
and $ ( \theta_{1} ,\theta_{2} ) \in [ - \pi ,\pi ] \times [ - \pi ,\pi ]$. 
Recalling (\ref{eq:sy}), it follows from (\ref{eq:ch2}), (i) in Lemma 5.1 and (iv) in Lemma 5.2 that 
the second term on the right-hand side of 
the last equality in (\ref{eq:b-b}) is bounded by a constant multiple of 
$  \{ \vert \theta_{1} \vert ^{ \hat{\delta}}+ 
\vert A \vert K^{(-1+ \delta )/4} + \vert C \vert 
( \vert \theta_{1} \vert ^{ \hat{\delta}} + K^{ \delta /4}) K^{-1/4} + 
\vert D \vert  K^{ (-2+ \delta )/4}  \} K^{-1/4} .$ 
Hence we have the desired estimate.

\par

\end{document}